\documentclass[11pt]{article}

\usepackage{graphicx}        
\usepackage{multicol}        
\usepackage{multirow}
\usepackage{amsmath,amssymb,amsthm}
\usepackage[british]{babel}
\usepackage[utf8]{inputenc}
\usepackage{xcolor}
\usepackage[linesnumbered,ruled,vlined,commentsnumbered]{algorithm2e}
\usepackage{cite}
\usepackage{url}
\usepackage{hyperref}
\usepackage{authblk}
\usepackage{caption}
\newtheorem{definition}{Definition}
\newcommand\rev[1]{\textcolor{black}{#1}}

\begin{document}

\title{Asymptotic-Preserving Neural Networks for hyperbolic systems with diffusive scaling}
\author{Giulia Bertaglia}
\affil{Department of Environmental and Prevention Sciences, University of Ferrara, via Luigi Borsari 46, 44121 Ferrara, Italy. \texttt{giulia.bertaglia@unife.it}}

\maketitle

\begin{abstract}
With the rapid advance of Machine Learning techniques and the deep \rev{increase} of availability of scientific data, data-driven approaches have started to become progressively popular across science, causing a fundamental shift in the scientific method after proving to be powerful tools with a direct impact in many areas of society. Nevertheless, when attempting to analyze dynamics of complex multiscale systems, the usage of standard Deep Neural Networks (DNNs) and even standard Physics-Informed Neural Networks (PINNs) may lead to incorrect inferences and predictions, due to the presence of small scales leading to reduced or simplified models in the system that have to be applied consistently during the learning process. In this Chapter, we will address these issues in light of recent results obtained in the development of \textit{Asymptotic-Preserving Neural Networks} (APNNs) for hyperbolic models with diffusive scaling. Several numerical tests show how APNNs provide considerably better results with respect to the different scales of the problem when compared with standard DNNs and PINNs, especially when analyzing scenarios in which only little and scattered information is available.
\end{abstract}

\tableofcontents

\section{Introduction}
\label{sec:intro}
The study of dynamics of complex systems described by multiscale partial differential equations (PDEs) has applications ranging from classical physics and engineering to biology, socio-economics and life sciences in general \cite{albi2022,pareschi2013,carrillo2010,bellomo2020,boscheri2021}. 
Despite the continuous progress achieved in the understanding of such systems, the modeling and prediction of the evolution of nonlinear multiscale phenomena using classical analytical or computational tools goes through several demanding challenges. Firstly, the numerical solution of a multiscale problem requires complex and elaborate computational codes and can introduce prohibitive costs, due to the well-known \textit{curse of dimensionality}. Secondly, especially when working in the social sciences, it is inevitable to come up against the difficulties associated with data scarcity and multiple sources of uncertainty \cite{albi2020,bertaglia2021,bertaglia2022,bertaglia2021a,bertaglia2021b}. Not to mention that solving real physical problems with missing or incomplete initial or boundary conditions is currently impractical with traditional approaches. 

With rapid advances in the field of Machine Learning (ML) and hugely increasing amounts of scientific data, data-driven approaches have started to become increasingly popular across science, causing a fundamental shift in the scientific method after proving to be tools with innumerable potentials and with direct impact in many areas of society \cite{higham2019,e2021,peng2021,baker2019,fabiani2021}.
One can then attempt to gain support from ML techniques even for the study of complex systems described by multiscale dynamics. When doing so, approaching the problem with the usage of standard ML methods, it is important to recall that purely data-driven models may fit observations very well (especially when a huge amount of data is available), but predictions may be physically inconsistent or unrealistic, with extrapolations leading to erroneous generalizations. Since, in general, we possess a-priori, albeit incomplete, knowledge of the physical laws governing the phenomenon under study, the key idea in the blossoming field of physics-informed ML is to include this type of prior scientific knowledge into the machine learning work-flow, ``teaching'' ML models not only to match observed data, but also to respect, in the extrapolations produced, the physical laws that we know govern the dynamics of interest. In fact, this approach aims to complement any (very likely) scarcity of available data with the knowledge of mathematical physical models, even in partially\rev{-}understood, uncertain and high-dimensionality contexts, guaranteeing robustness, more accurate and, above all, physically consistent predictions \cite{lou2021}. 

A recent example reflecting this new learning philosophy is represented by Physics-Informed Neural Networks (PINNs), first introduced in \cite{raissi2019}. PINNs constitute a new class of deep neural networks (DNNs) that are trained to solve supervised learning tasks respecting physical laws described through linear or nonlinear ordinary differential equations (ODEs) or PDEs. The physical knowledge of the phenomenon under study, in addition to being incorporated into the learning process of the DNN directly through data embodying the underlying physics of the phenomenon of interest (observational bias), is introduced through an appropriate choice of the loss function that the DNN has to minimize, which forces the training phase of the neural network to converge towards solutions that adhere to the physics (learning bias)~\cite{karniadakis2021,coutinho2022,yu2021,moseley2020}.

However, the use of a standard formulation of PINNs in the context of multiscale problems can still lead to incorrect inferences and predictions \cite{bertaglia2022a,jin2021}. This is mainly due to the presence of small scales leading to reduced or simplified models in the system that have to be applied consistently during the learning process. In these cases, a standard PINN formulation only allows an accurate description of the process at the leading order, thus losing accuracy in asymptotic limit regimes. A remedy to this problem, as recently proposed in \cite{bertaglia2022a,jin2021}, is to modify the loss function to include asymptotic-preserving (AP) properties during the training process. The realization of such an AP-loss function will therefore depend on the particular problem under consideration and be based on an appropriate asymptotic analysis of the model. 

In this Chapter, we will address these issues in the light of recent results obtained in the development of \emph{asymptotic-preserving neural networks} (APNNs) for hyperbolic models with diffusive scaling. More precisely, the presentation is organized as follows: 
Section \ref{sec:hyp_model} is devoted to the overview of hyperbolic systems with relaxation terms leading to an asymptotic parabolic limit; 
in Section \ref{sec:NN-PINN} we review the basic concepts of DNNs and PINNs; 
in Section \ref{sec:APNN} the construction of APNNs with AP loss functions for multiscale hyperbolic systems of interest is discussed, giving an analytical demonstration of the AP property and its importance in the context of neural networks; 
in Section \ref{sec:numericsGT} the different behavior of a standard DNN, a standard PINN and an APNN are analyzed and compared considering the Goldstein-Taylor model as prototype system of equations, evaluated for different regimes; 
finally, an application of APNNs in the context of multiscale hyperbolic transport models for the study of the propagation of infectious diseases is presented in Section \ref{sec:SIR} and conclusive remarks are given in Section \ref{sec:conclusion}.

\section{Hyperbolic systems with diffusive scaling}
\label{sec:hyp_model}
The prototype hyperbolic system with relaxation that we will consider is a discrete-velocity kinetic model describing the space-time evolution of two particle densities, $f^+(x,t)$ and $f^-(x,t)$, at time $t>0$, traveling in a 1D domain, $x\in\mathcal{D}\subseteq\nobreak\mathbb{R}$, with opposite velocities $\pm 1$, respectively, which is known as Goldstein-Taylor model \cite{pareschi2013,jin2000,jin2015}:
\begin{eqnarray}
\frac{\partial f^+}{\partial t} + \frac{1}{\epsilon} \frac{\partial f^+}{\partial x} &=& -\frac{\sigma}{2\epsilon^2} \left(f^+ - f^-\right)\,, \label{eq.GT1}\\
\frac{\partial f^-}{\partial t} - \frac{1}{\epsilon} \frac{\partial f^-}{\partial x} &=& \frac{\sigma}{2\epsilon^2} \left(f^+ - f^-\right)\,.
\label{eq.GT2}
\end{eqnarray}
Here $\sigma$ identifies a \rev{possible} scattering coefficient and $\epsilon$ represents the scaling factor. Indeed, in this model, particles can collide and, at the same time, assume the opposite velocity, randomly, with a probability that depends on the scaling factor $\epsilon$, which is strictly related to the \textit{mean free path} of particles: the more $\epsilon$ tends to zero, the smaller the mean free path of the particles and the more collisions will occur. In this context, the role of the operator which describes the interactions with the background is to push the solution towards a universal steady state (at exponential-in-time rate proportional to the collision frequency). The typical definition used in kinetic theory for this type of operator is indeed \textit{relaxation operator} \cite{pareschi2013}.

If we define the total particles density $\rho=f^++f^-$ and the re-scaled flux~$j=\nobreak\epsilon^{-1}(f^+-f^-)$, we can rewrite system \eqref{eq.GT1}--\eqref{eq.GT2} obtaining an equivalent model expressed in terms of macroscopic variables:
\begin{eqnarray}
\frac{\partial \rho}{\partial t} +  \frac{\partial j}{\partial x} &=& 0\,,	\label{eq.hypsyst.1}\\ 
\frac{\partial j}{\partial t} + \frac{1}{\epsilon^2} \frac{\partial \rho}{\partial x} &=& -\frac{\sigma}{\epsilon^2}j \,,
\label{eq.hypsyst.2}
\end{eqnarray}
With this formulation, we can analyze the behavior of the solution as $\epsilon \to 0$, i.e., in the \textit{zero-relaxation limit}. In this asymptotic limit, the second equation relaxes to the local equilibrium 
\begin{equation}
j= - \frac1{\sigma}\frac{\partial \rho}{\partial x}\,,
\label{eq.j_diff}
\end{equation} 
and substituting into the first equation gives the following parabolic, diffusive limit that recalls the standard heat equation:
\begin{equation}
\frac{\partial \rho}{\partial t}  = \frac{\partial}{\partial x}\left[\frac1{\sigma}\frac{\partial \rho}{\partial x}\right] \,.
\label{eq.heat}
\end{equation}
It is therefore clear that, depending on the scaling factor $\epsilon$, system \eqref{eq.GT1}--\eqref{eq.GT2} describes different distinct physical propagation phenomena that range from that of advective transport to that of parabolic diffusion \cite{albi2019,boscarino2017}, associated with spatio-temporal scales which differ by several orders of magnitude.

In classical kinetic theory, the space-time scaling just discussed is related to one specific hydrodynamical limit of the Boltzmann equation. In particular, when the dissipative effects become non negligible, from the Boltzmann equation we recover the incompressible Navier--Stokes scaling. We refer the reader to \cite{cercignani1994} for further details and the mathematical theory behind the hydrodynamical limits of the Boltzmann equation and to \cite{lions1997} for theoretical results on the diffusion limit of a system like \eqref{eq.hypsyst.1}--\eqref{eq.hypsyst.2}.

\section{Review of Deep Neural Networks and Physics-Informed Neural Networks}
\label{sec:NN-PINN}
In this Section, we provide a brief summary of the general framework of DNNs and their extension to PINNs \cite{raissi2019,karniadakis2021}. After this overview, the relevant concepts of  APNNs for the problems of interest will be discussed in Section \ref{sec:APNN}.

\subsection{Deep Neural Networks (DNNs)}
\label{sec:NN}
\begin{figure}[t!]
\centering
    \includegraphics[width=0.55\textwidth]{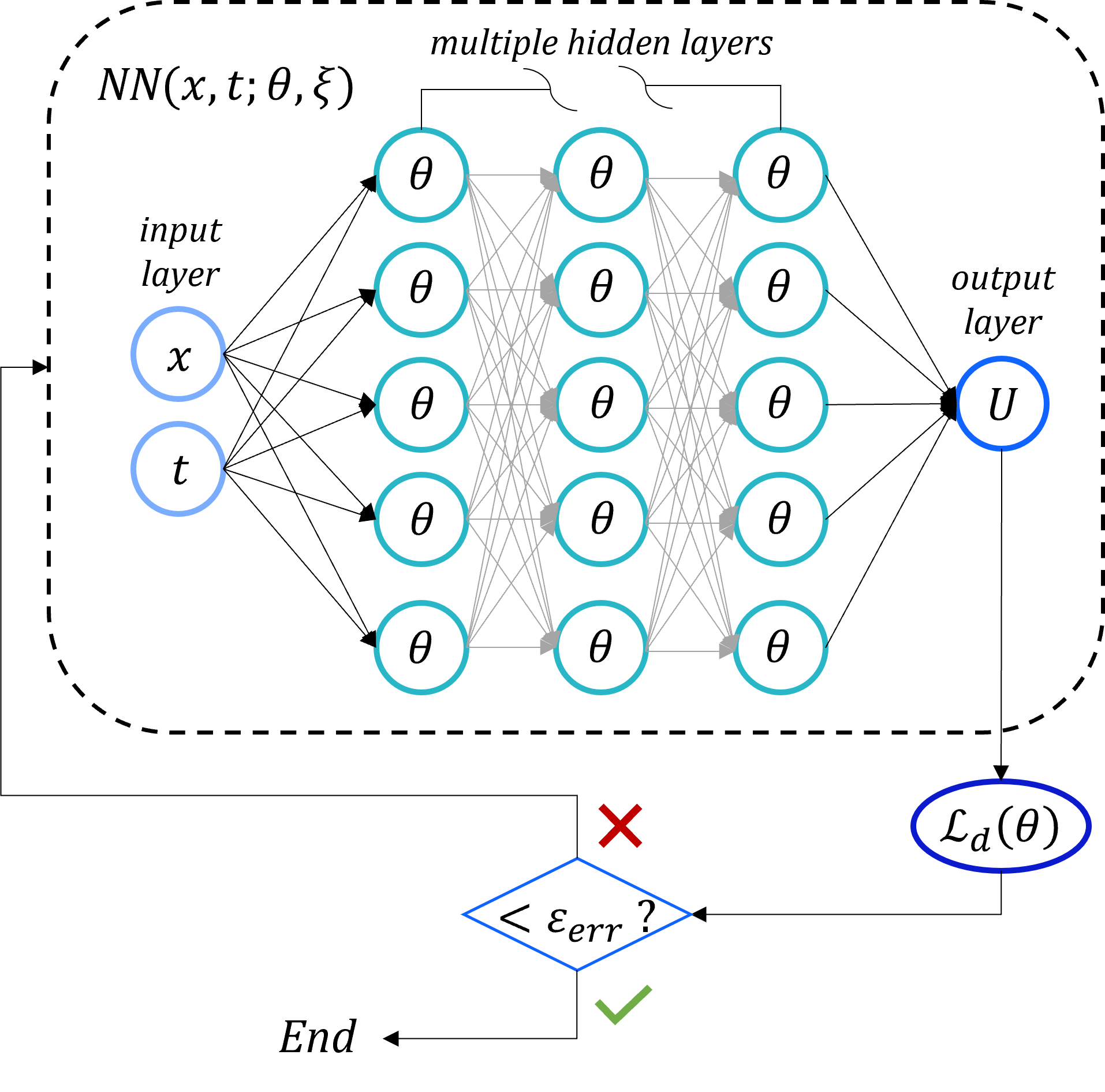}
    \caption{NN schematic work-flow. Given the location $(x,t)$ of the available training data as input layer, to find the optimal values for $\theta$ in the hidden layers, the neural network is trained by minimizing the loss function $\mathcal{L}_d(\theta)$ that generally consists in the mean squared error between the NN's predictions (output layer) and the training points. Iterations end when the error evaluated at the training points is less than a set threshold $\varepsilon_{\mathit{err}}$.}
    \label{fig:NN_scheme}
\end{figure}

Let us assume we want to evaluate the dynamic of a variable $U\in \Omega \subset \mathbb{R}^{d}\times\mathbb{R}$ through a data-driven approach, recurring to a classical DNN, for example an $L+1$ layered feed-forward neural network (FNN), which consists of an input layer, an output layer, and $L-1$ hidden layers. Given the location of an available dataset of $N_d$ observations $z\in\mathbb{R}^d, d=2,$ as input, the FNN can be used to output a prediction of the value $U \approx U_{\mathit{NN}}(z; \theta)$, parameterized by network parameters $\theta$, as shown in Figure \ref{fig:NN_scheme} for a spatio-temporal dynamic in 1D space. To this aim, we can define the FNN as follows:
\begin{eqnarray*}
z^1 &=& W^1 z + b^1\,,\\
z^l &=& \sigma \circ (W^l z^{l-1} + b^l),\quad l=2,\ldots,L-1\,,\\
U_{\mathit{NN}}(z;\theta) = z^L&=& W^L z^{L-1} + b^L\,,
\end{eqnarray*}
where $W^{l}\in\mathbb{R}^{m_{l}\times m_{l-1}}$ are the weights, $b^l\in\mathbb{R}^{m_{l}}$ the bias, $m_l$ is the width of the $l$-th hidden layer with $m_1 = d_{in} = d$ the input dimension and $m_L = d_{out}$ the output dimension,  $\sigma$ is a scalar activation function (such as ReLU \cite{goodfellow2016}), and ``$\circ$'' denotes entry-wise operation. Thus, we denote the set of network parameters~$\theta=(W^1,b^1,\ldots,W^L,b^L)$.

To learn the model, the network’s free parameters $\theta$ are tuned to find the optimal values at each iteration of the neural network through a supervised learning process, so that the DNN’s predictions closely match the available experimental data. This is usually done by minimizing a loss function $\mathcal{L}_d(\theta)$ (also called cost or risk function) that consists in the mean squared error (MSE) between the DNN's predictions and the training points, which for a spatio-temporal dynamic in 1D space reads:
\begin{equation}
\label{eq:general-nn-loss}
\mathcal{L}_d(\theta) = \frac{1}{N_d} \sum_{i=1}^{N_d} \left| U_{\mathit{NN}}(x_d^i, t_d^i;\theta) - U(x_d^i, t_d^i) \right|^2\,,
\end{equation}
where ${U}_d^i , i=1,\ldots,N_d$ is the available measurement dataset and $\{(x_d^i, t_d^i)\}_{i=1}^{N_d} \subset \Omega$ are the coordinates of these training points.

This problem is typically cast as a stochastic optimization problem and generally solved using a stochastic gradient descent (SGD) algorithm such as the \rev{Adam} advanced optimizer \cite{Adam}.
After finding the optimal set of parameter values $\theta^*$ by minimizing the loss function \eqref{eq:general-nn-loss}, i.e.,
\begin{equation*}
\theta^* = \mathrm{argmin}\, \mathcal{L}_d(\theta),
\end{equation*}
the neural network surrogate $U_{\mathit{NN}}(x,t; \theta^*)$ can be evaluated at any point in the domain to get the solution of the problem.
 
Therefore, the design of a DNN can be summarized in the following three main steps \cite{E2020}: 
\begin{itemize}
\item the choice of the neural network structure,
\item the choice of the loss function,
\item the choice of the method to minimize the loss over the parameter space.
\end{itemize}

Notice that, in practice, the performance of the neural network is estimated on a finite set of points which is unrelated to any data used to train the model evaluating a \emph{test error} (or \emph{validation error}); whereas the error in the loss function, which is used for training purposes, is called the \emph{training error}.

\subsection{Physics-Informed Neural Networks (PINNs)}
\label{sec:PINN}
As already discussed in Section \ref{sec:intro}, using a purely data-driven approach like the one just presented can have significant weaknesses. Indeed, whilst a standard DNN might accurately model the physical process within the vicinity of the experimental data, it will fail to generalize away from the training points. This is because the neural network cannot actually \emph{learn} the physical dynamic of the phenomenon \cite{moseley2020}, and recurring to standard DNNs we are essentially forgetting about our existing scientific knowledge. 

\begin{figure}[t!]
\centering
    \includegraphics[width=\textwidth]{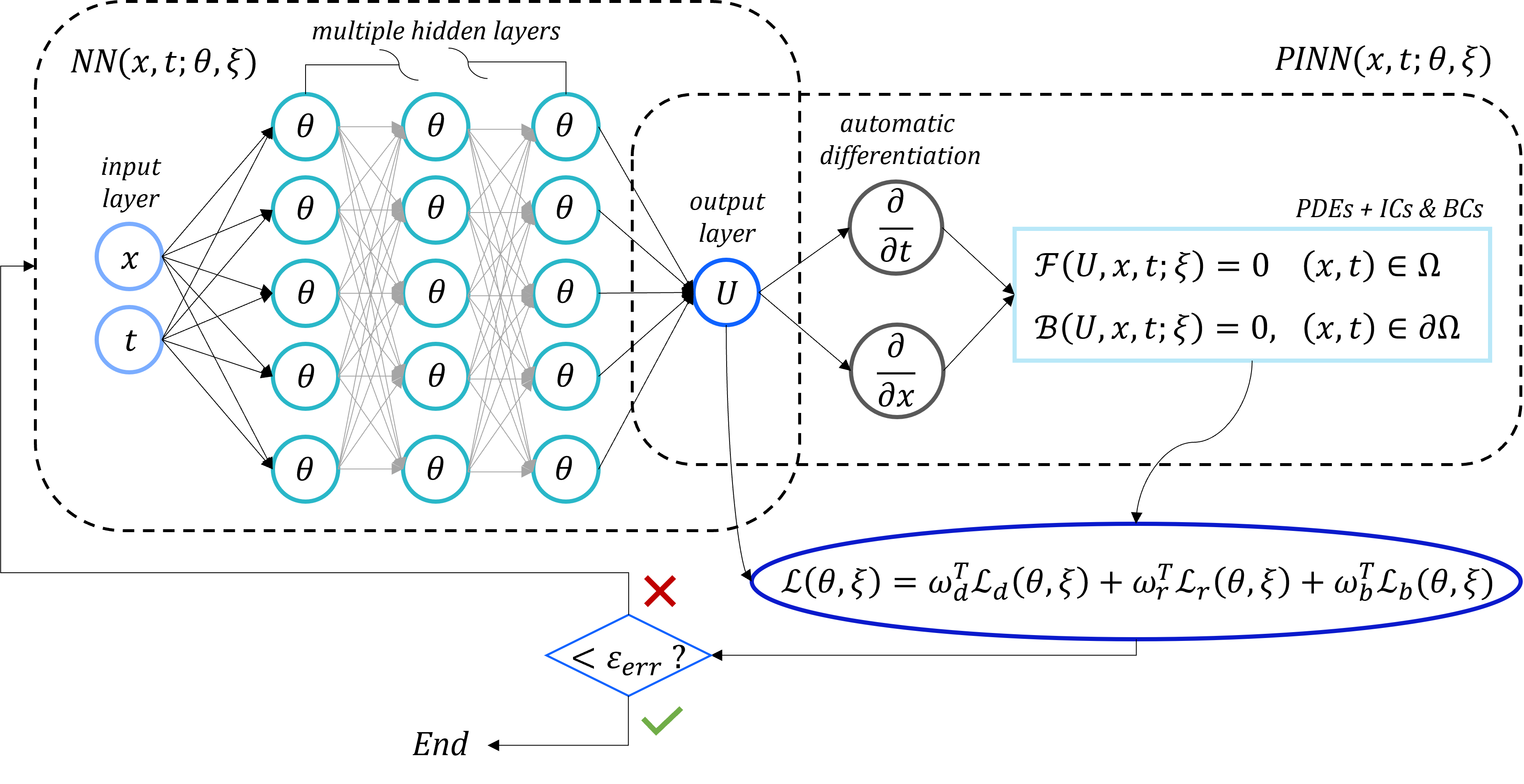}
    \caption{PINN schematic work-flow. The NN architecture is integrated with the physical knowledge of the dynamic of interest through the inclusion of the differential operator (in this case a PDE) and the enforcement of initial and boundary conditions (and eventually conservation properties), when known, becoming a PINN. To this aim, the loss function encloses three different residual terms, related to the mismatch with respect to available data $\mathcal{L}_d(\theta)$, known physical law $\mathcal{L}_r(\theta)$ and boundary/initial conditions $\mathcal{L}_b(\theta)$, respectively.}
    \label{fig:PINN_scheme}
\end{figure}
To improve the efficiency and efficacy of the learning process, PINNs take advantage of the prior physical knowledge in our \rev{possession} including it into the neural network work-flow, providing a powerful tool that can address both the inverse and the forward problem of data-driven learning of the dynamic and solving the problem of interest, respectively \cite{raissi2019,karniadakis2021,kharazmi2021}. 

More specifically, let us first consider the context of inverse problems assuming that, in addition to having access to experimental data, we know that the dynamic investigated is governed by a differential operator $\mathcal{F}$ in the spatio-temporal domain~$\Omega\subset\mathbb{R}^{d}\times\mathbb{R}$ which depends on parameters $\xi$ related to the physics whose values are unknown. Furthermore, the dynamic is subject to a general operator $\mathcal{B}$ that prescribes arbitrary initial and boundary conditions of the system in $\partial \Omega$, which could be known at specific locations. This generic governing system will read as follows:
\begin{eqnarray}
\label{eq:general-pde-system}
    \mathcal{F}(U, x, t; \xi) &=& 0, \quad  (x, t) \in \Omega, \\ \label{eq:general-pde-system_BC}
    \mathcal{B}(U, x, t; \xi) &=& 0, \quad (x, t) \in \partial \Omega.
\end{eqnarray}

To transform a standard NN into a PINN, it is sufficient to require the neural network to satisfy not only the agreement with experimental data\rev{,} but also the physics we know of the system, acting directly \rev{on} the loss function by adding two more terms: 
\begin{equation}
\label{eq:general-pinn-loss}
\mathcal{L}(\theta,\xi) = w_d^{\rev{T}} \mathcal{L}_{d}(\theta,\xi) + w_b^{\rev{T}} \mathcal{L}_{b}(\theta,\xi) +  w_r^{\rev{T}} \mathcal{L}_{r}(\theta,\xi),
\end{equation}
with
\begin{equation}
\mathcal{L}_{r}(\theta,\xi) = \frac{1}{N_r} \sum_{n=1}^{N_r} \left| \mathcal{F}(U_{\mathit{NN,r}}^n,x_r^n, t_r^n;\theta,\xi) \right|^2\,,
\end{equation}
\begin{equation}
\mathcal{L}_{b}(\theta,\xi) = \frac{1}{N_b} \sum_{k=1}^{N_b} \left| \mathcal{B}(U_{\mathit{NN,b}}^k,x_b^k, t_b^k;\theta,\xi) \right|^2\,,
\end{equation}
where $\{(x_r^n, t_r^n)\}_{n=1}^{N_r} \subset \Omega$ and $\{(x_b^k, t_b^k)\}_{k=1}^{N_b} \subset \partial \Omega$ are scattered points within the domain and on the boundary, respectively, called \emph{residual points}, at which PINN is asked to satisfy the known physics; while ${U}_{\mathit{NN,r}}^n~=~{U}_{\mathit{NN}}(x_r^n, t_r^n)$ and ${U}_{\mathit{NN,b}}^k~=~{U}_{\mathit{NN}}(x_b^k, t_b^k)$.
Here, the additional residual terms, $\mathcal{L}_{r}$ and $\mathcal{L}_{b}$, quantify the discrepancy of the neural network surrogate $U_{\mathit{NN}}$ with respect to the underlying differential operator and its initial or boundary conditions in \eqref{eq:general-pde-system_BC}, respectively. Finally, $w_{r}$, $w_{b}$, $w_{d}$ are the weights associated to each contribution. \rev{Notice that if one aims at penalizing a specific loss term with respect to the others (e.g., $\mathcal{L}_r$, in the context of a scalar $\mathcal{F}$ and $\mathcal{B}$) it is either possible to set a weight smaller than one in front of it (e.g., $\omega_r = 0.1$, $\omega_d=\omega_b=1$) or to use weights higher than one in front of the rest of the loss terms in the loss function $\mathcal L$ (e.g., $\omega_d=\omega_b=10$, $\omega_r = 1$).} 

In this additional process, gradients of the network’s output with respect to its input are computed at each residual point typically recurring to automatic differentiation~\cite{baydin2017}, and the residual of the underlying differential equation is evaluated using these gradients.

In the context of inverse problems, the unknown physical parameters $\xi$ are treated as learnable parameters, and we search for those that best describe the observed data. As a result, the training process involves optimizing $\theta$ and $\xi$ jointly:
\begin{equation*}
    (\theta^*, \xi^*) = \mathrm{argmin}\,  \mathcal{L}(\theta,\xi).
\end{equation*}
For a schematic representation of the discussed PINN work-flow the reader can refer to Figure \ref{fig:PINN_scheme}.

In the context of forward problems, the structure of the network is almost the same as in the inverse problem, except that no matching with data in the learning phase is normally considered and we assume to know $\xi$. In essence, we ask \rev{the} PINN what can be said about the unknown hidden state $U$ of the system, given the fixed parameters~$\xi$ of the model.  Thus, the optimization task will generally read:
\begin{equation*}
\theta^* = \mathrm{argmin}\, \left[w_r^{\rev{T}}\mathcal{L}_{r}(\theta) + w_b^{\rev{T}}\mathcal{L}_{b}(\theta)\right] .
\end{equation*}

The physics-informed neural networks thus constructed are able to predict the solution also away from the experimental data points by behaving much better than a standard neural network, as we will explicitly see in Section \ref{NNvsPINN}.

\section{Asymptotic-Preserving Neural Networks}
\label{sec:APNN}
\begin{figure}[b!]
\centering
    \includegraphics[width=0.6\textwidth]{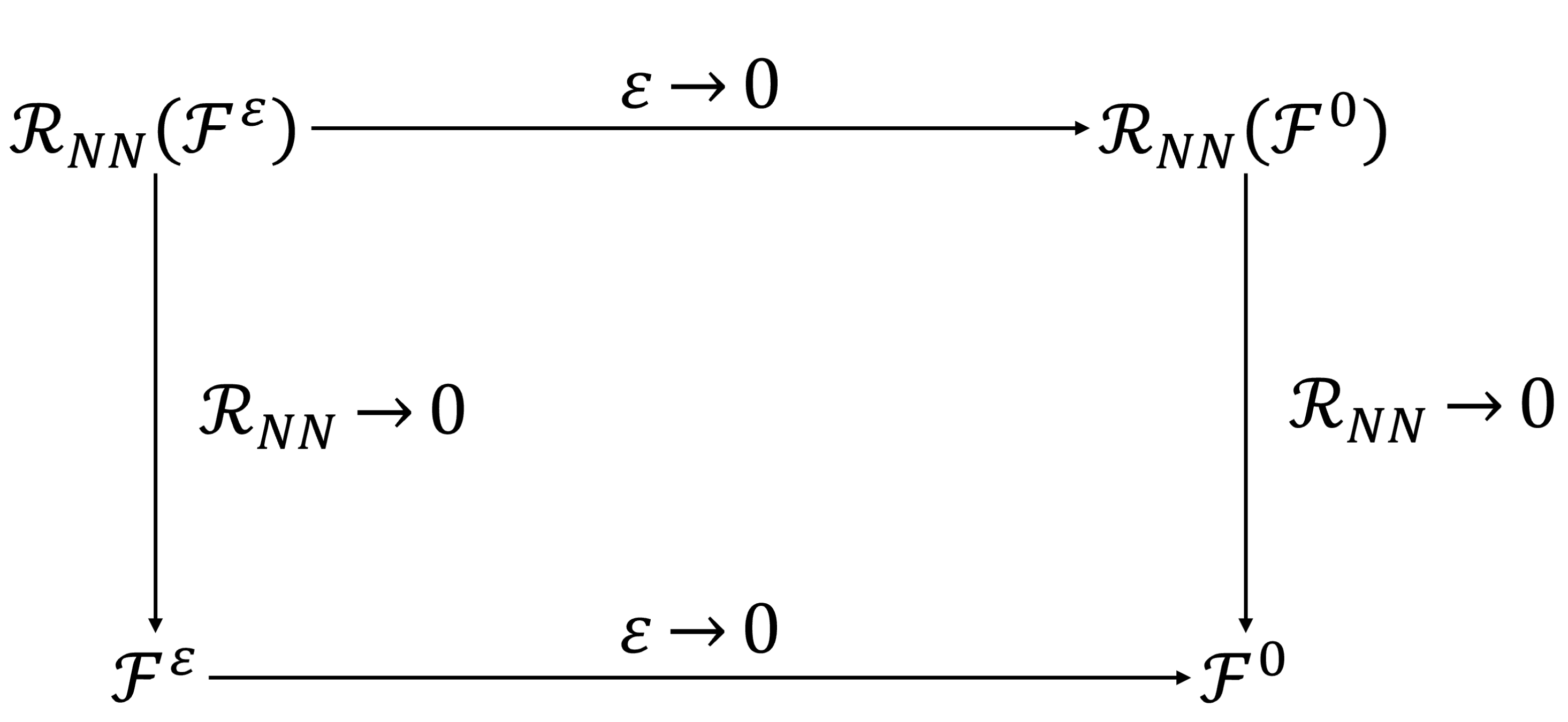}
    \caption{AP diagram for physics-informed neural networks. $\mathcal{F}^{\epsilon\mkern-2mu}$ is the multiscale hyperbolic model that depends on the scaling parameter $\epsilon$, while $\mathcal{F}^0$ is the corresponding formulation in the diffusive limit, which does not depend anymore on $\epsilon$. The solution of the system $\mathcal{F}^\epsilon$ is approximated by the neural network through the imposition of the residual term $\mathcal{R}_{\mathit{NN}}\left( \mathcal{F}^{\epsilon\mkern-2mu}\right) = \mathcal{R}_{\mathit{NN}}^{\epsilon\mkern-2mu}$. The asymptotic limit of $\mathcal{R}_{\mathit{NN}}\left( \mathcal{F}^{\epsilon\mkern-2mu}\right)$ as $\epsilon \to 0$ is denoted with $\mathcal{R}_{\mathit{NN}}\left( \mathcal{F}^0\right) = \mathcal{R}_{\mathit{NN}}^0$. The neural network is called AP if $\mathcal{R}_{\mathit{NN}}\left( \mathcal{F}^{\epsilon\mkern-2mu}\right)$ represents a good approximation of $\mathcal{F}^0$.}
    \label{fig:AP_scheme}
\end{figure}

Since we aim at analyzing multiscale hyperbolic dynamics regardless of the propagation scaling, in order to obtain physically-based predictions, it is important that the PINN is able to preserve also the correct equilibrium solution in the diffusive regime of the equations, such \rev{as} that of system \eqref{eq.hypsyst.1}--\eqref{eq.hypsyst.2} with limit \eqref{eq.heat}.
The neural networks satisfying this requirement are called {\it Asymptotic-Preserving Neural Networks}~(APNNs), and have been recently introduced in \cite{jin2021,jin2022} precisely to efficiently solve multiscale kinetic problems with scaling parameters that can have several orders of magnitude of difference. In essence, an APNN is said to be such if the neural network (already PINN) benefits from the \emph{asymptotic preservation} (AP) property, by analogy with that of classical numerical schemes. In this analogy, the loss function is viewed as a numerical approximation of the original equation and has to benefit from the AP property. Let us recall that the so-called AP schemes constitute a class of numerical methods that aim to preserve the correct asymptotic behavior of the system without any loss of efficiency due to time step restrictions related to the small scales \cite{albi2020,boscarino2017,dimarco2014}.
Therefore, the definition of an APNN (reported in \cite{jin2021} for the case of multiscale kinetic models with continuous velocity fields) can be generalized as follows (see Fig. \ref{fig:AP_scheme}).

\begin{definition}
Assume the solution of the problem is parametrized by a PINN trained by using an optimization method to minimize a loss function which includes a residual term $\mathcal{R}_{\mathit{NN}}$ enforcing the physics of the phenomenon $\mathcal{F}^\epsilon\mkern-2mu$, which depends on a scaling parameter $\epsilon$. Then we say it is an \textit{Asymptotic-Preserving Neural Network} (APNN) if, as the physical scaling parameter of the multiscale model tends to zero (i.e., $\epsilon \to 0$), the loss function of the full model-constraint converges to the loss function of the corresponding reduced order model $\mathcal{F}^0$. 
\end{definition}

\subsection{APNN for the Goldstein-Taylor model}
\label{AP-GT}
To illustrate the relevance of the AP property in the construction of a neural network for hyperbolic systems with diffusive scaling, let us carry on a detailed example by directly considering the Goldstein-Taylor model \eqref{eq.GT1}--\eqref{eq.GT2}. 

We assume ${f}_{\mathit{NN}}^\pm(x,t;\theta,\sigma)$ to be the outputs of a PINN given inputs $x$ and $t$ and trainable parameters $\theta$ and $\sigma$ (so we assume the scattering coefficient to be unknown), which approximate the solution of our system: ${f}^\pm(x,t;\sigma) \approx {f}_{\mathit{NN}}^\pm(x,t;\theta,\sigma)$. Then, we define the PDEs residuals by multiplying both members of the equations \rev{by} the \rev{square of the} scaling parameter $\epsilon\rev{^2}$ to allow the usage of the model also when $\epsilon \to 0$,
\begin{equation}
    \label{eq:pinn-nonap-residue1}
    \mathcal{R}^{\epsilon\mkern-2mu, \pm}_{\mathit{NN}} = \epsilon^2\frac{\partial f^\pm_{\mathit{NN}}}{\partial t} \pm \epsilon \frac{\partial f^\pm_{\mathit{NN}}}{\partial x} - \frac{\sigma}{2} \left(f^\mp_{\mathit{NN}}-f^{\pm}_{\mathit{NN}}\right),
\end{equation}
and incorporate them into the loss function term ${\mathcal L}_r(\theta,\sigma)$ of the neural network in~\eqref{eq:general-pinn-loss} by taking the mean squared error to obtain a standard PINN:
\begin{equation}
\label{eq:PINNloss}
\omega_r^T {\mathcal L}_r(\theta,\sigma) =\frac{\omega_r^+}{N_r} \sum_{n=1}^{N_r} \left| \mathcal{R}_{\mathit{NN}}^{\epsilon\mkern-2mu,+} (x_r^n, t_r^n;\theta,\sigma) \right|^2 +
  \frac{\omega_r^-}{N_r} \sum_{n=1}^{N_r} \left| \mathcal{R}_{\mathit{NN}}^{\epsilon\mkern-2mu,-} (x_r^n, t_r^n;\theta,\sigma) \right|^2\,.
\end{equation}

It is easy to observe that, with this construction, the standard PINN residual \eqref{eq:pinn-nonap-residue1} is not consistent with the reduced order model of the system, namely the diffusion limit \eqref{eq.heat}. In fact, $\mathcal{R}^{\epsilon,\pm}_{\mathit{NN}}$ in the limit $\epsilon\to 0$ reduce to 
\begin{equation*}
\mathcal{R}^{0,\pm}_{\mathit{NN}} = - \frac{\sigma}{2} \left(f^\mp_{\mathit{NN}}-f^{\pm}_{\mathit{NN}}\right)\,,
\end{equation*}
which basically corresponds to \rev{forcing} $f^+(x,t)=f^-(x,t)$ and does not suffice to achieve the correct diffusive behavior. Therefore, the thus constructed PINN does not satisfy the sought AP property.

In contrast, using the macroscopic formulation of the model \eqref{eq.hypsyst.1}--\eqref{eq.hypsyst.2} results sufficient to construct an APNN, with outputs ${\rho}_{\mathit{NN}}(x,t;\theta,\sigma)$ and ${j}_{\mathit{NN}}(x,t;\theta,\sigma)$. Incorporating now in the loss function the mean squared error of the PDEs residuals of the macroscopic formulation,
\begin{equation}
    \label{eq:pinn-ap-residue_GT}
    \mathcal{R}^{\epsilon\mkern-2mu,\rho}_{\mathit{NN}} =  \frac{\partial \rho_{\mathit{NN}}}{\partial t} +  \frac{\partial j_{\mathit{NN}}}{\partial x}\,,\qquad
    \mathcal{R}^{\epsilon\mkern-2mu, j}_{\mathit{NN}} =  \epsilon^2 \frac{\partial j_{\mathit{NN}}}{\partial t} +  \frac{\partial \rho_{\mathit{NN}}}{\partial x}+\sigma j_{\mathit{NN}}\,,
\end{equation}
the loss term of the physics residual will read
\begin{equation}
\label{eq:APNNloss}
\omega_r^T {\mathcal L}_r(\theta,\sigma) =\frac{\omega_r^\rho}{N_r} \sum_{n=1}^{N_r} \left| \mathcal{R}_{\mathit{NN}}^{\epsilon\mkern-2mu,\rho} (x_r^n, t_r^n;\theta,\sigma) \right|^2 +
  \frac{\omega_r^j}{N_r} \sum_{n=1}^{N_r} \left| \mathcal{R}_{\mathit{NN}}^{\epsilon\mkern-2mu,j} (x_r^n, t_r^n;\theta,\sigma) \right|^2\,.
\end{equation}
In the limit $\epsilon \to 0$, we obtain
\begin{equation*}
    \mathcal{R}^{0,\rho}_{\mathit{NN}} =  \frac{\partial \rho_{\mathit{NN}}}{\partial t} +  \frac{\partial j_{\mathit{NN}}}{\partial x}\,,\qquad
    \mathcal{R}^{0,j}_{\mathit{NN}} =  \frac{\partial \rho_{\mathit{NN}}}{\partial x}+\sigma j_{\mathit{NN}}\,,
\end{equation*}
which analytically confirms that we are now dealing with a loss function that remains consistent also with the residual of the limiting diffusive model \eqref{eq.heat}.

\section{Application to the Goldstein-Taylor model}
\label{sec:numericsGT}
In this Section, we will analyze and compare the behavior of a standard DNN with respect to a standard PINN and of a standard PINN and an APNN when trying to solve the Goldstein-Taylor model in a non-multiscale and a multiscale configuration, respectively. 
\subsection{Standard DNN vs Standard PINN in hyperbolic regime}
\label{NNvsPINN}
\begin{figure}[t!]
    \centering
     \includegraphics[width=0.45\textwidth]{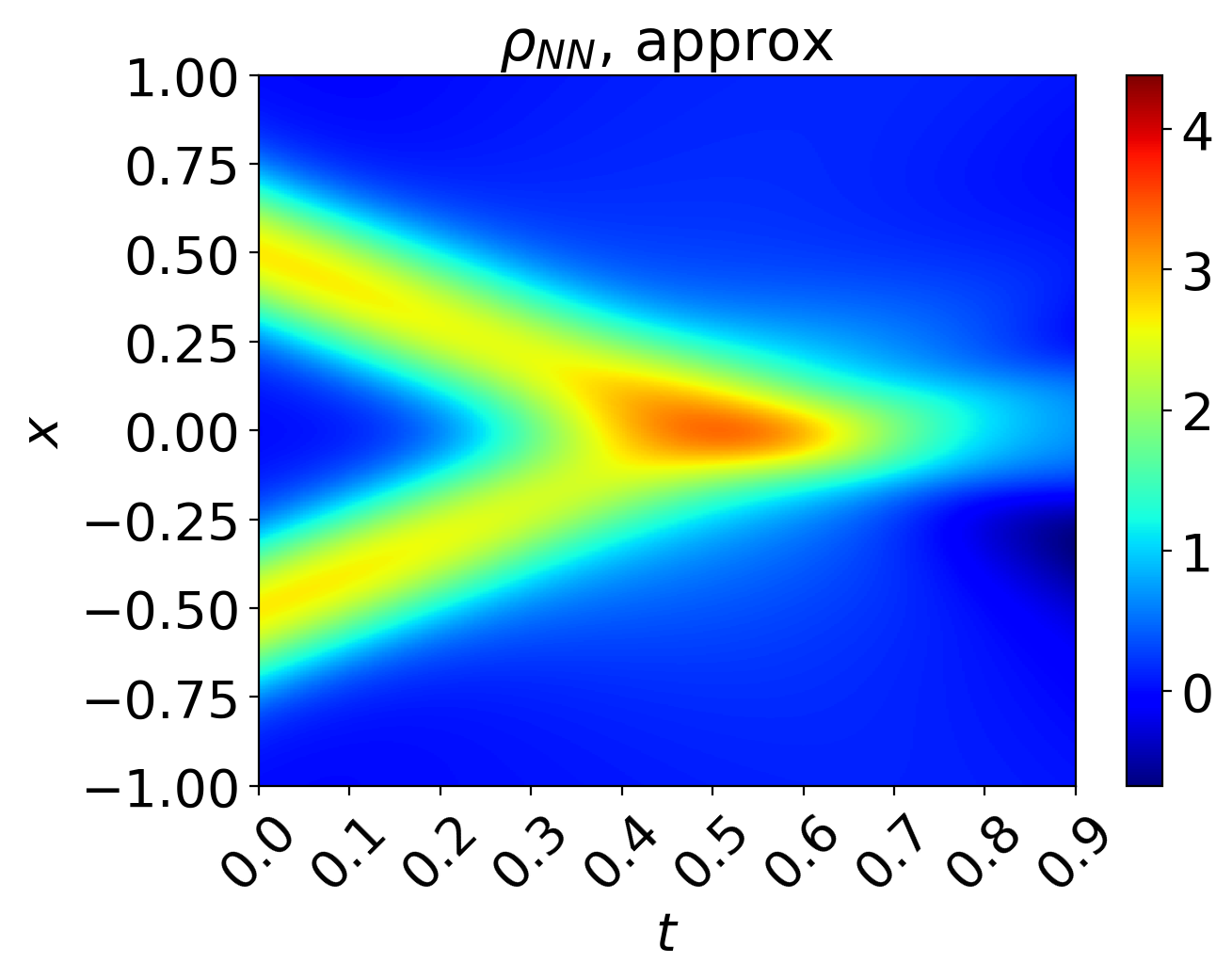}
     \includegraphics[width=0.45\textwidth]{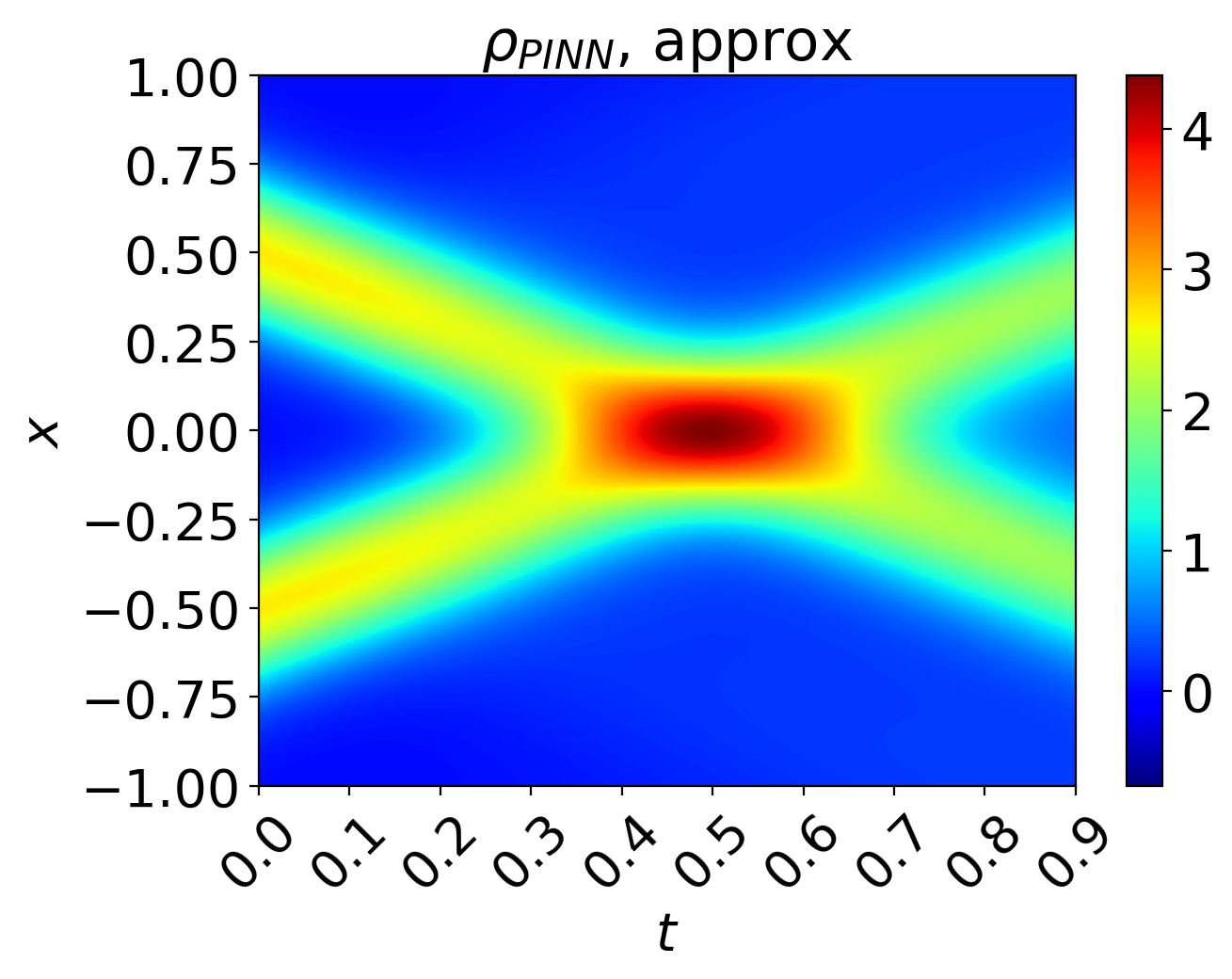}
     \includegraphics[width=0.42\textwidth]{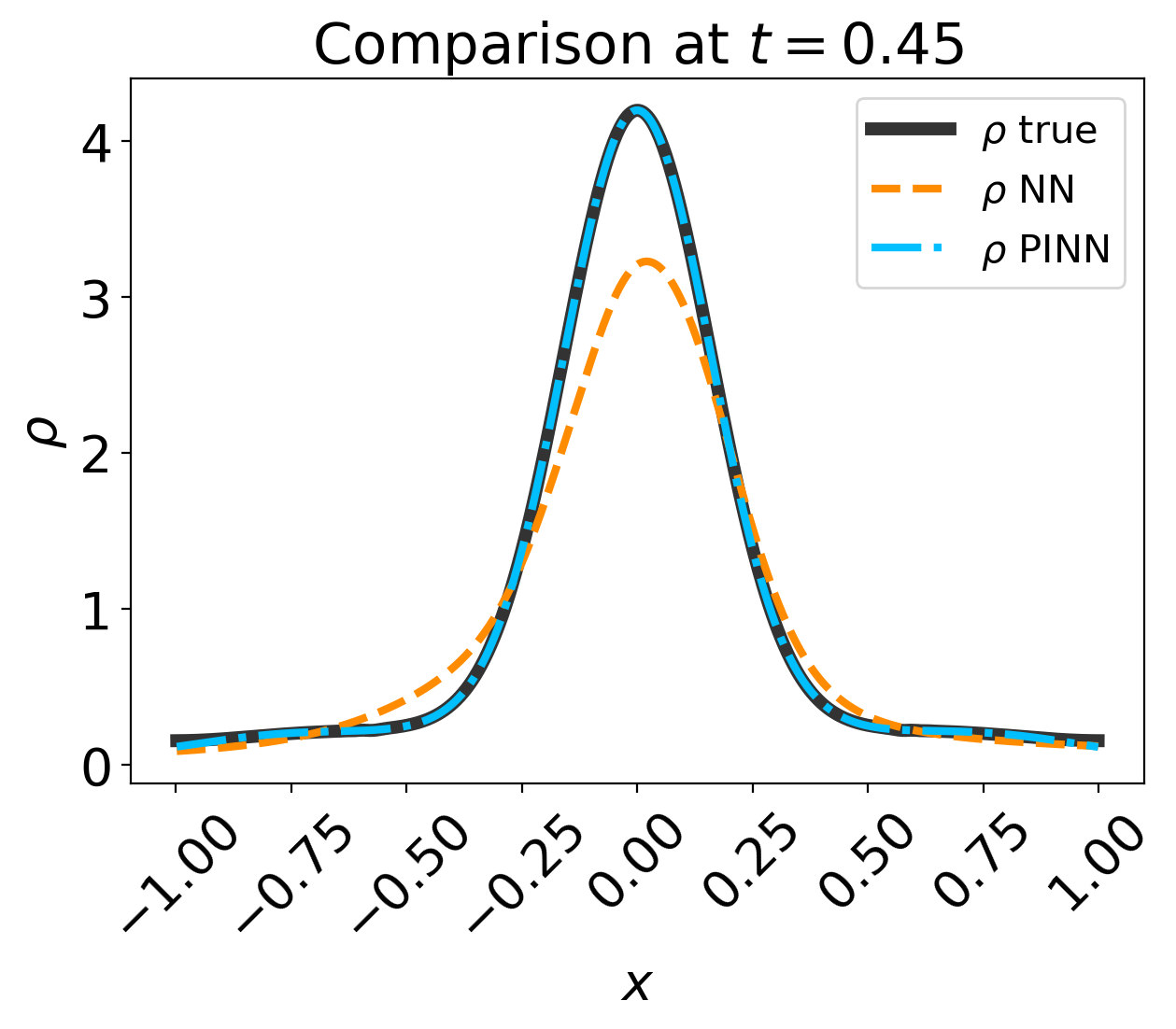}
     \hspace{0.01\textwidth}
     \includegraphics[width=0.45\textwidth]{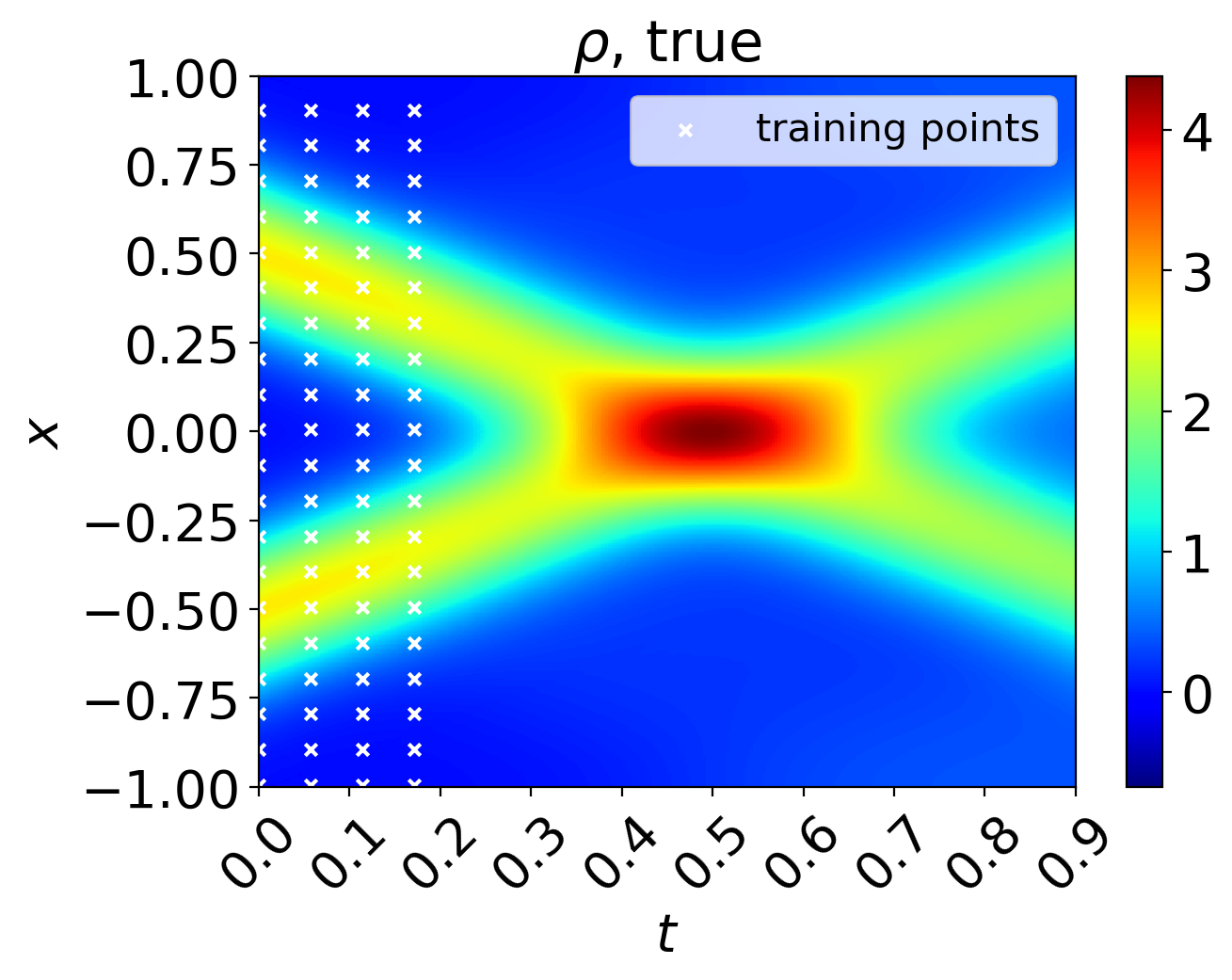}
    \caption{Forward problem for the Goldstein-Taylor model with $\epsilon=1.0$. Solution of the forward problem obtained with a standard DNN (top, left) or with a standard PINN (top, right) with respect to the ground truth (bottom, right) in terms of total density $\rho$; comparison of the DNN and PINN solution at $t=0.45$ (bottom, left). White crosses in the bottom right plot identify the 80 training data points used for both the neural networks.}
    \label{fig:rho_hypGT}
\end{figure}

To first analyze the different behavior of a standard DNN with respect to a PINN, let us consider a numerical test for the Goldstein-Taylor model presented in Section~\ref{sec:hyp_model} without the presence of the scaling parameter, hence fixing $\epsilon=1.0$. Therefore, we build up a DNN that is trained only with observed data points taken from a synthetic solution obtained solving system \eqref{eq.GT1}--\eqref{eq.GT2} with a second-order AP-IMEX Runge-Kutta finite volume method \cite{bertaglia2021b,bertaglia2020b} that is considered as ground truth (see Figure \ref{fig:rho_hypGT} \rev{bottom,} right). On the other hand, we consider also a PINN having in the loss function the additional physics loss term as in \eqref{eq:PINNloss}. Moreover, in the PINN\rev{'s} loss function we enforce the positivity of the densities by adding the term
\begin{eqnarray*}
\mathcal{L}_p(\theta) &=&\frac{1}{N_r} \sum_{n=1}^{N_r} \left| \mathrm{abs}\left(f_{\mathit{NN}}^{+} (x_r^n, t_r^n;\theta)\right) -  f_{\mathit{NN}}^{+} (x_r^n, t_r^n;\theta)\right|^2 \\&+&  \frac{1}{N_r} \sum_{n=1}^{N_r} \left| \mathrm{abs}\left(f_{\mathit{NN}}^{-} (x_r^n, t_r^n;\theta)\right) -  f_{\mathit{NN}}^{-} (x_r^n, t_r^n;\theta)\right|^2 \,.
\end{eqnarray*}

In the domain $\mathcal{D}=[-1,1]$, we consider two initial Gaussian distributions of densities $f^\pm$ as
$$
f^\pm(x,0) = \frac{1}{s\sqrt{2\pi}} e^{-\frac{(x\pm0.5)^2}{2s^2}}
$$
with $s = 0.15$ and scattering coefficient $\sigma=1$, which is assumed to be known. We consider periodic boundary conditions and run the simulation until $t_{\mathit{end}}=0.9$.

In both the NNs, we consider a feed-forward network with depth $3$ (i.e., number of layers, including output layer but excluding input layer) and width $32$ (i.e., number of nodes in a layer). We use $\tanh$ as activation function and fix a learning rate LR=$10^{-3}$. Finally, the Adam method \cite{Adam} is employed for the optimization process and the derivatives in the physics loss function are computed applying automatic differentiation \cite{baydin2017}.
The weights in the PINN loss function \eqref{eq:general-pinn-loss} are chosen in order to penalize the physics residuals. Hence, we set all of them equal to 1, except for $\omega_d=10$.
The spatio-temporal approximation of the total density $\rho$ obtained with the standard DNN trained with $N_d=80$ observation points of densities $f^+$ and $f^-$ is compared with the one obtained using a standard PINN, which is trained not only \rev{on} the same dataset but also enforcing the PDE system structure and the density positivity on $N_r=3600$ residual points uniformly distributed in the entire domain (without enforcing boundary conditions) in Figure \ref{fig:rho_hypGT}. From this Figure, it is clearly visible that the standard DNN is not able to reproduce a correct solution of the problem when trained with limited data and does not fulfil to predict credible trends far from the dataset availability. In contrast, the PINN is able to fit the solution not only in the vicinity of the training points\rev{,} but also outside of them, giving reliable predictions.

\subsection{Standard PINN vs APNN in diffusive regime}
To emphasize also numerically the importance of choosing the correct formulation to fulfill the AP property and correctly approximate dynamics of multiscale systems even in diffusive regimes, we consider a test similar to the one discussed in the previous Section for the Goldstein-Taylor model, but choosing $\epsilon = 10^{-5}$ and evaluating the solution until $t_{\mathit{end}}=0.02$. 
In the domain $\mathcal{D}=[-1,1]$, centered in $x_0=0$, we consider a problem with an initial Gaussian distribution of the total density $\rho$ \rev{that is the same obtained in the previous test 
$$
\rho(x,0) = \frac{1}{s\sqrt{2\pi}} \left(e^{-\frac{(x-0.5)^2}{2s^2}} + e^{-\frac{(x+0.5)^2}{2s^2}}\right)\,,
$$
but evaluating the corresponding equilibrium flux $j$ with \eqref{eq.j_diff}, which results}
$$
\rev{ j(x,0) = - \frac{1}{\sigma s^3\sqrt{2\pi}} \left( (x-0.5)e^{-\frac{(x-0.5)^2}{2s^2}} + (x+0.5)e^{-\frac{(x+0.5)^2}{2s^2}} \right)\,, }
$$
again with $s = 0.15$, $\sigma=1$ (assumed to be known), and periodic boundary conditions.

We solve the problem first applying the standard PINN, with the same setting previously discussed, and then, for comparison, using the APNN formulation, with residuals \eqref{eq:pinn-ap-residue_GT} in the AP loss function and enforcing the positivity of the total density:
\begin{equation*}
\mathcal{L}_p(\theta) =\frac{1}{N_r} \sum_{n=1}^{N_r} \left| \mathrm{abs}\left(\rho_{\mathit{NN}} (x_r^n, t_r^n;\theta)\right) -  \rho_{\mathit{NN}} (x_r^n, t_r^n;\theta)\right|^2  \,.
\end{equation*}
Notice that the standard PINN shares the same setting with the APNN, except for the different residual terms enforced for the physics, as discussed in Section \ref{AP-GT}. In the APNN, indeed, we consider the correctly scaled system \eqref{eq:APNNloss} instead of \eqref {eq:PINNloss}.

\begin{figure}[!t]
    \centering
     \includegraphics[width=0.45\textwidth]{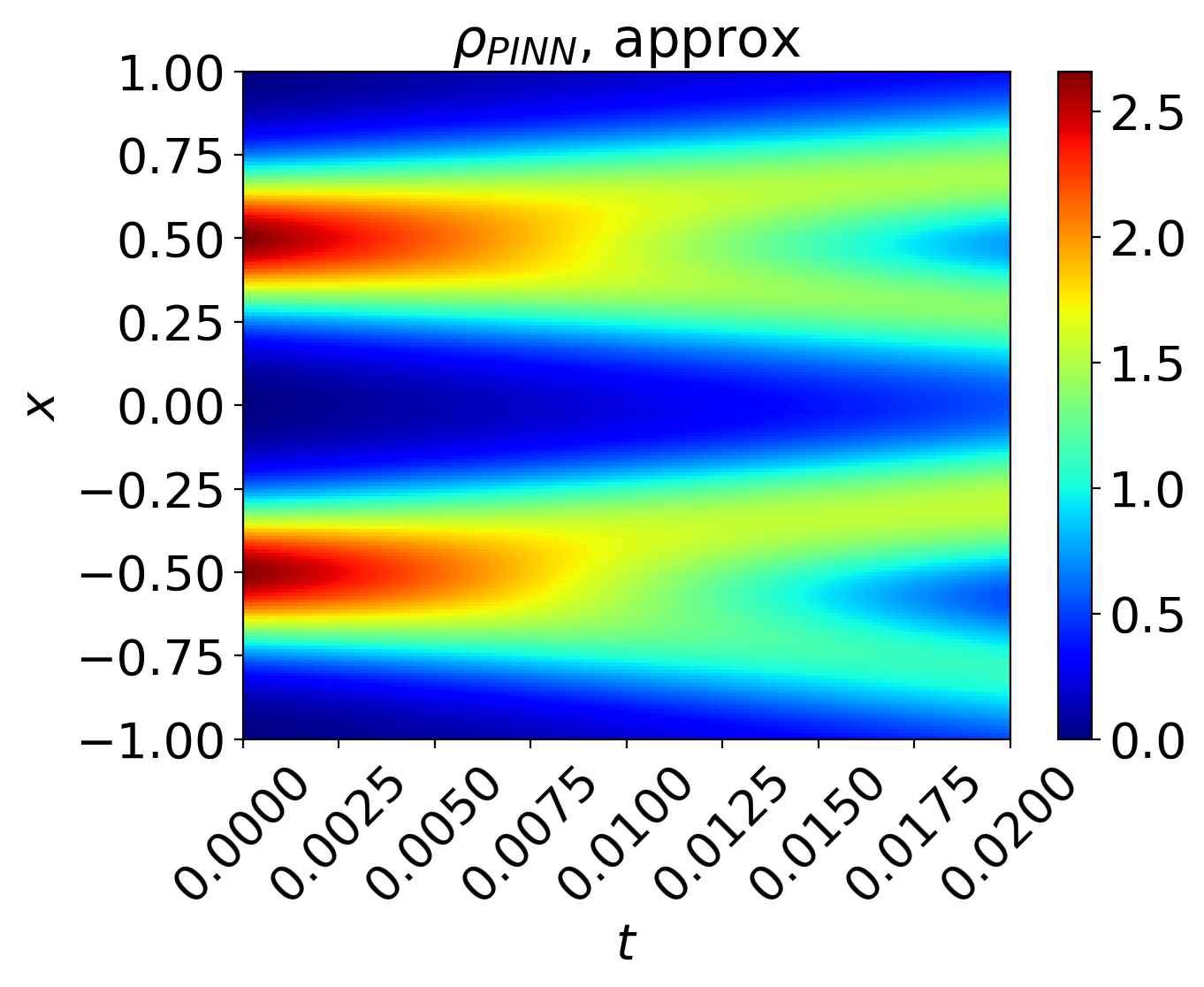}
     \includegraphics[width=0.45\textwidth]{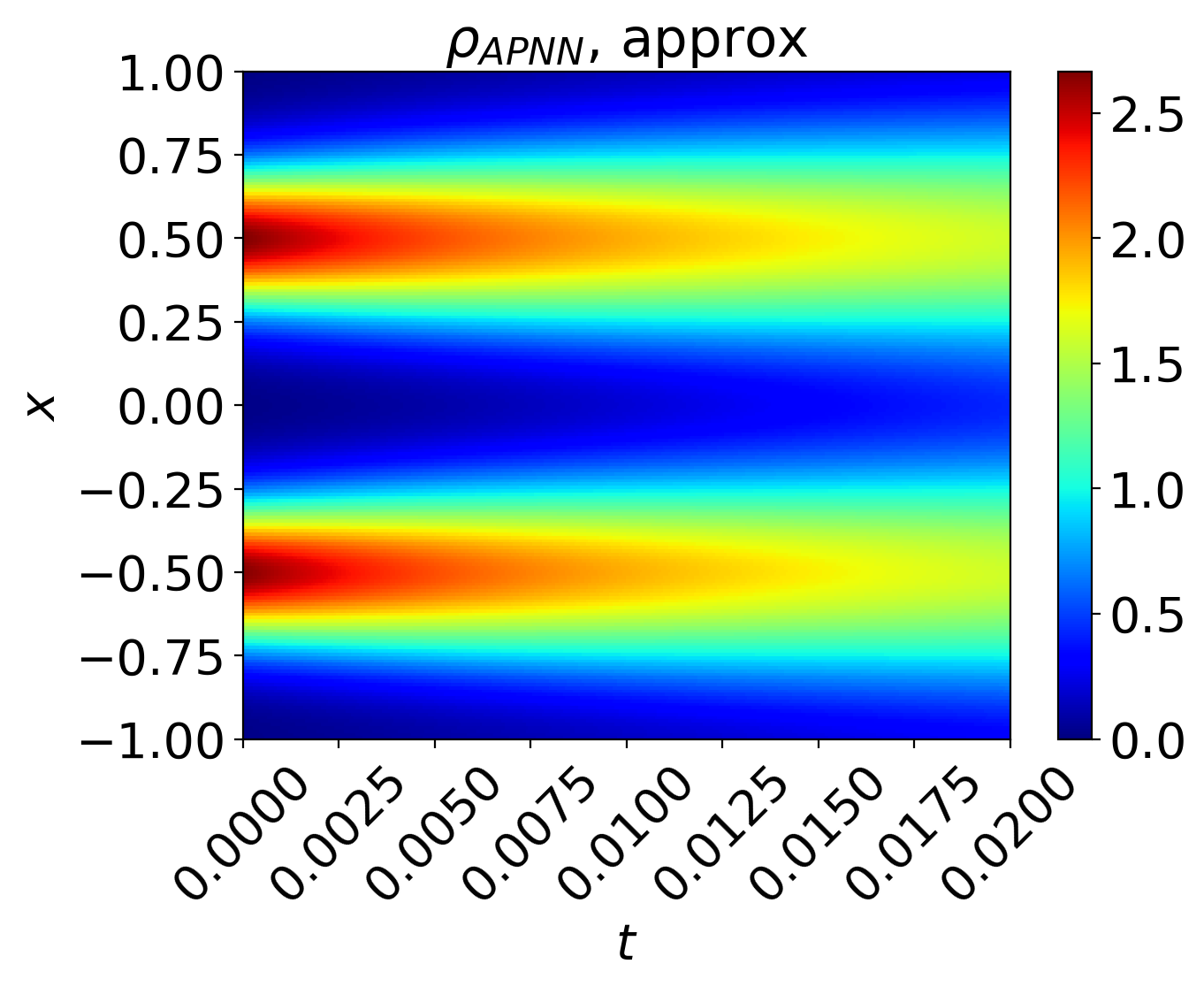}
      \includegraphics[width=0.45\textwidth]{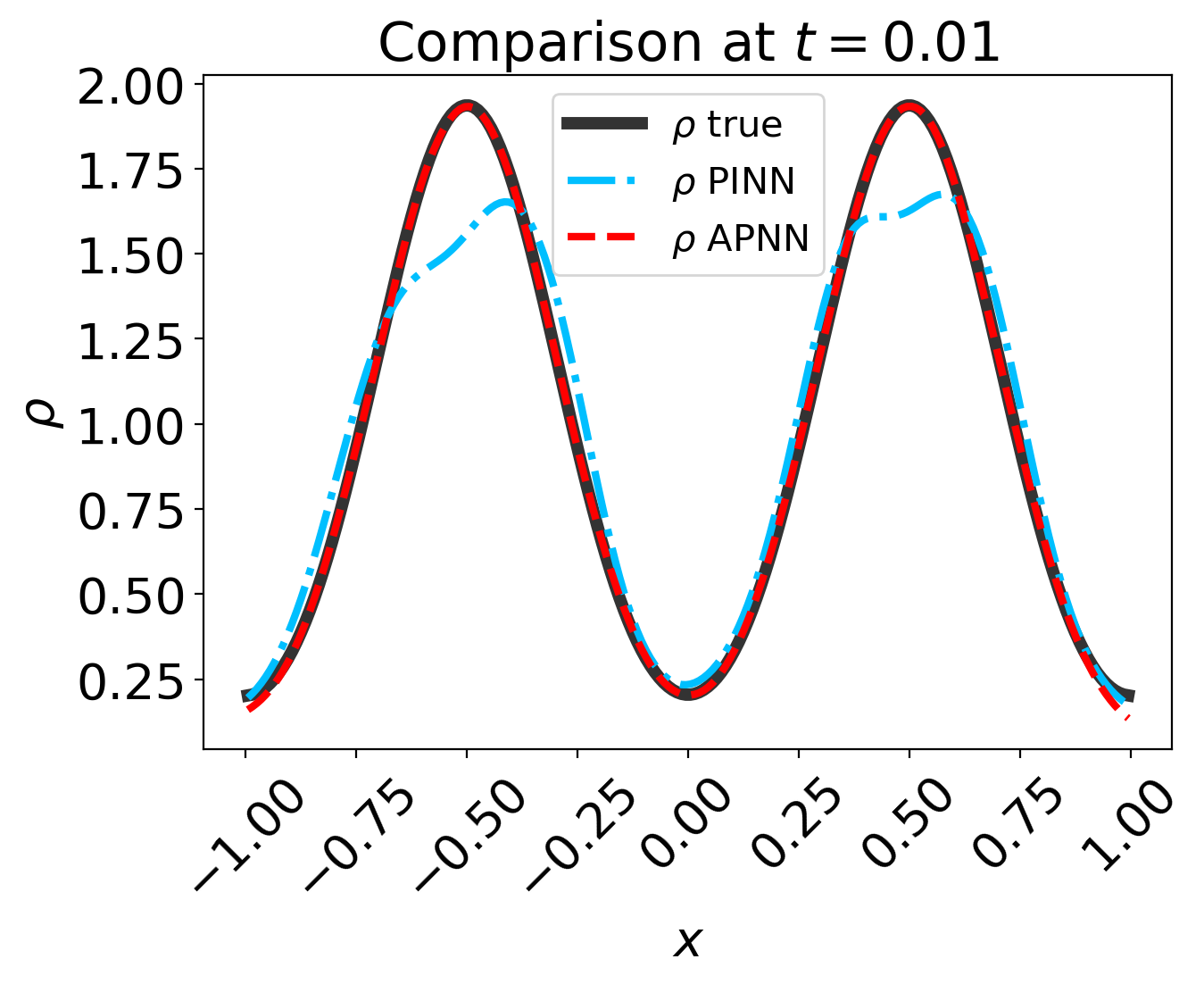}
     \includegraphics[width=0.45\textwidth]{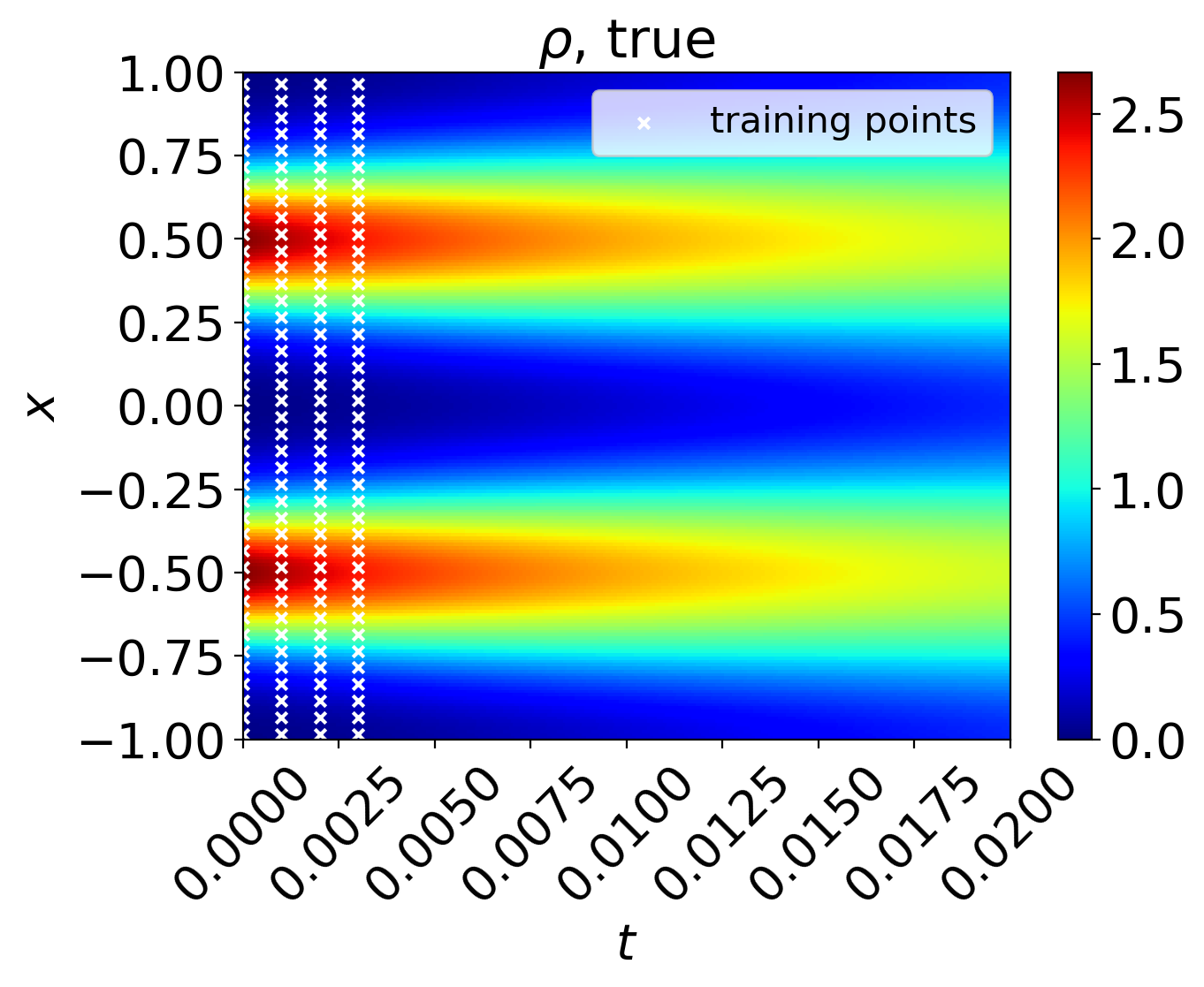}
    \caption{Forward problem for the Goldstein-Taylor model in diffusive regime ($\epsilon=10^{-5}$). Solution of the forward problem obtained with a standard PINN (top, left) or with an APNN (top, right) with respect to the ground truth (bottom, right) in terms of total density $\rho$; comparison of the PINN and APNN solution at $t=0.01$ (bottom, left). White crosses in the right plot identify the 80 sparse training data points used to train both the neural networks.}
    \label{fig:rho_diffGT}
\end{figure}
Results obtained with the two neural networks, trained with $N_d=160$ data points and enforcing the two different physics residuals in $N_r=2500$ uniformly distributed points, are presented in Figure \ref{fig:rho_diffGT}. We observe that the standard PINN is not able to recover the correct solution in this regime, especially with respect to forecast scenarios far from the dataset, because the physics imposed (and learned by the PINN) is not consistent with the relaxation limit of the system. In contrast, the APNN approximates the solution very well even at temporal points where there is no nearby data availability.

\paragraph{Inverse Problem}
To analyze also the behavior of standard PINN and APNN when dealing with the resolution of inverse problems, we design an additional test for the Goldstein-Taylor model, for which initial conditions are
\begin{equation*}
\rho(x,0) = 6 + 3\cos(3\pi x)\,, \qquad j(x,0) = 9\pi \sigma^{-1}\sin(3\pi x)\,,
\end{equation*}
with $\sigma=4$, model parameter that is assumed to be unknown and need to be estimated by the neural network. 
We consider again periodic boundary conditions and we fix $\epsilon=10^{-4}$ to simulate the diffusive, parabolic regime of the model until $t_{\mathit{end}}=0.1$.

For both standard PINN and APNN formulations, we train the network model on measurements composed of $N_d=24000$ equally spaced samples (taken from the synthetic solution) in the domain $(x, t)\in [-1, 1]\times [0, 0.1]$, from which 20\%~(4800) points are randomly selected for validation purpose. For the APNN model, to account for an even poorer partially available dataset, we consider measurements only for the density $\rho$, hence assuming to have no information on the flux $j$, whereas for the standard PINN formulation we employ data samples for both the densities~$f^+$ and~$f^-$ (therefore, in the latter case we assume we have more information on the system). Finally, $N_r=24000$ residual points are employed to enforce the physics, with the same 20\% data split for the validation set. 
\begin{figure}[!t]
    \centering
    \includegraphics[width=.48\textwidth]{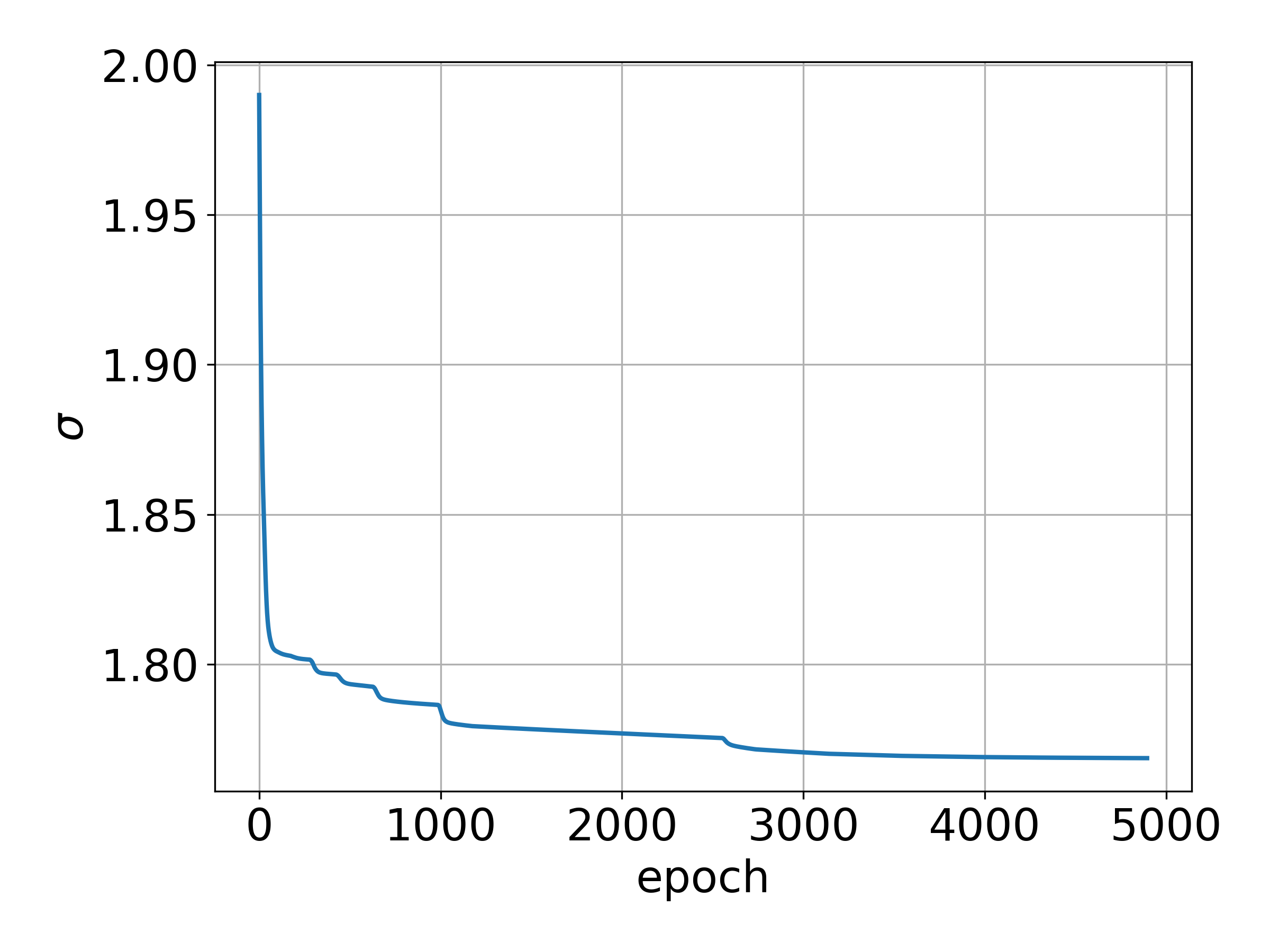}
    \includegraphics[width=.48\textwidth]{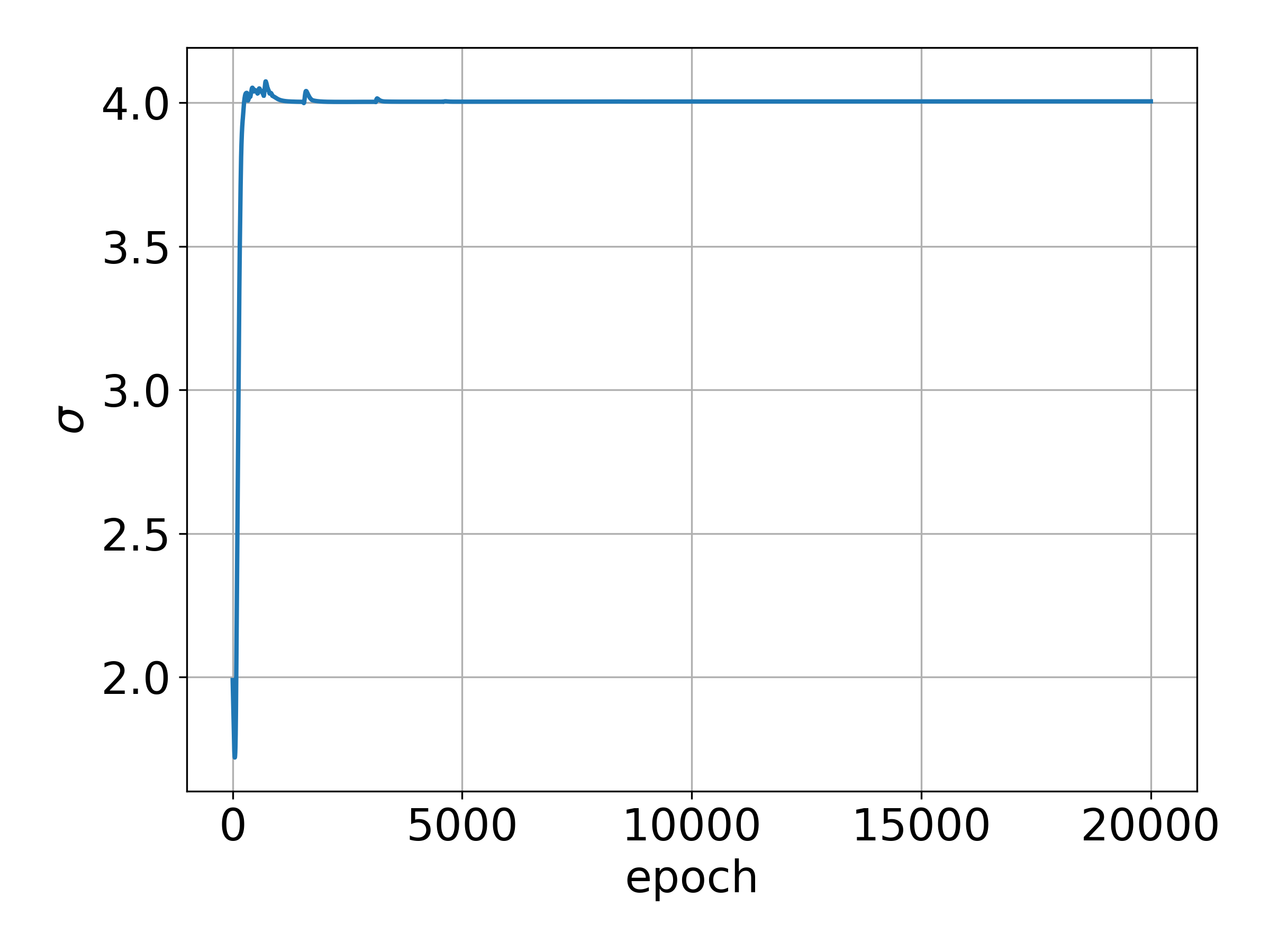}  
    \caption{Inverse problem for the Goldstein-Taylor model in diffusive regime ($\epsilon=10^{-4}$). Convergence of the target parameter $\sigma=4$ with respect to epochs using the standard PINN (left) and the APNN~(right). The correct inference of the parameter is obtained only when recurring to the APNN.}
    \label{fig:test1_parabolic_sigma}
\end{figure}

We show the convergence of the target parameter $\sigma$ in Figure \ref{fig:test1_parabolic_sigma} for both the neural networks. Despite the partial availability of information on the dynamic (of which flux data are completely missing), a very fast convergence of the APNN towards the expected value of $\sigma$ is observed, obtaining a final relative error $\mathcal{O}(10^{-3})$. On the other hand, it can be observed that the standard PINN, although trained with a full dataset of $f^\pm$, fails to retrieve the correct value of the scattering coefficient $\sigma$ (early stopping the iterations almost at epoch 5000 to avoid further training). 

\section{Application to epidemic dynamics}
\label{sec:SIR}
A particularly interesting area in which the use of the APNN framework discussed here can play a key role is that of the study of epidemic dynamics. In this context, also driven by the COVID-19 pandemic that has hit the world, a number of mathematical models have recently been proposed that require the estimation of several parameters from data to provide predictive scenarios and to test their reliability. Examples are~\cite{albi2022,bertaglia2020b,albi2020,boscheri2021,buonomo2020,gatto2020, giordano2021,marca2022,scarabel2021,guglielmi2021}.
Here we focus on a new class of epidemic models defined by multiscale PDEs capable of describing both hyperbolic-type phenomena, characteristic of epidemic propagation over long distances and main lines of communication between cities, and parabolic-type phenomena, in which classical diffusion prevails at the urban level \cite{boscheri2021,bertaglia2021,bertaglia2021b}.

\subsection{The multiscale hyperbolic SIR model}
Let us consider for simplicity the space dependent epidemiological modeling in the case of a classic SIR compartmental dynamic. Thus, we start subdividing the population in susceptible $S$ (individuals who may be infected by the disease), infectious $I$ (individuals who may transmit the disease) and removed $R$ (individuals healed and immune or died due to the disease), following the pioneering work of Kermack and McKendrick \cite{kermack1927}. We assume that no subjects have prior immunity and we neglect the vital dynamic represented by births and deaths due to the time scale considered. 
By analogy with discrete-velocity kinetic theory \cite{pareschi2013}, and in particular with the Goldstein-Taylor model  \eqref{eq.GT1}--\eqref{eq.GT2}, we consider now that individuals are moving in a 1D domain $\mathcal{D} \subseteq \mathbb{R}$ in two opposite directions, with velocities $\pm \lambda_{\mathit{S,I,R}} = \pm \lambda_{\mathit{S,I,R}}(x)$, distinguished for each epidemic compartment. Notice that the characteristic velocities reflect the heterogeneity of geographical areas, and, therefore, are chosen dependent on the spatial location $x \in \mathcal{D}$. We can describe the space-time dynamic of this population for $t>0$ through the following two-velocity SIR epidemic transport model~\cite{bertaglia2021,bertaglia2020b}:
\begin{eqnarray}
	\frac{\partial S^{\pm}}{\partial t} + \lambda_S \frac{\partial S^{\pm}}{\partial x} &=& - \beta S^{\pm}I \mp \frac{1}{2\tau_S}\left(S^+ - S^-\right)\,,	\label{eq.SIR_kinetic_diagS}\\ 
	\frac{\partial I^{\pm}}{\partial t} + \lambda_I \frac{\partial I^{\pm}}{\partial x} &=&  \beta S^{\pm}I -\gamma I^{\pm} \mp \frac{1}{2\tau_I}\left(I^+ - I^-\right)\,, \label{eq.SIR_kinetic_diagI}\\ 
	\frac{\partial R^{\pm}}{\partial t} + \lambda_R \frac{\partial R^{\pm}}{\partial x} &=& \gamma I^{\pm} \mp \frac{1}{2\tau_R}\left(R^+ - R^-\right)\,,				
\label{eq.SIR_kinetic_diagR}
\end{eqnarray}
with the total densities of each compartment, $S(x,t)$, $I(x,t)$, and $R(x,t)$, given by
\begin{equation*}
 	S = S^+ + S^- \,, \quad 
 	I = I^+ + I^- \,, \quad 
 	R = R^+ + R^- \,,	
\end{equation*}
and 
$S(t)+I(t)+R(t)=P$, total population size, which remains constant over time.

The transport dynamic of the population is governed by the scaling parameters~$\lambda_{\mathit{S,I,R}}$ as well as the relaxation times $\tau_{\mathit{S,I,R}}=\tau_{\mathit{S,I,R}}(x)$. The quantity $\gamma=\gamma(x,t)$ is the recovery rate of infected, which corresponds to the inverse of the infectious period. 
The transmission of the infection is defined by an incidence function~$\beta S I$ modeling the transmission of the disease \cite{hethcote2000,capasso1978,kermack1927}, where the transmission rate~$\beta=\beta(x,t)$ characterizes the average number of contacts per person per time, multiplied by the probability of disease transmission in a contact between a susceptible and an infectious subject. 

When investigating real epidemic scenarios, the  above-mentioned  parameters are, in general, unknown. While the recovery rate might be fixed based on clinical data, the transmission rate must always be estimated through a delicate calibration process in order to match observations. It is also well-known that this process is highly heterogeneous, which makes the inverse problem even more challenging \cite{bertaglia2022,dimarco2021}.

If we look at system \eqref{eq.SIR_kinetic_diagS}--\eqref{eq.SIR_kinetic_diagR}, it is worth observing that, for each compartment, if we exclude the presence of the epidemic source term, we are essentially prescribing a dynamic that directly recalls that presented for the Goldstein-Taylor model \eqref{eq.GT1}--\eqref{eq.GT2} in Section \ref{sec:hyp_model}. 
Furthermore, defining now the flux for each compartment as
\begin{equation*}
 	J_S = \lambda_S \left(S^+ - S^-\right) \,, \quad 
 	J_I = \lambda_I \left(I^+ - I^-\right) \,, \quad 
 	J_R = \lambda_R \left(R^+ - R^-\right) \,,	  	
\end{equation*}
in analogy with the Goldstein-Taylor model, it is possible to derive the following macroscopic formulation of system \eqref{eq.SIR_kinetic_diagS}--\eqref{eq.SIR_kinetic_diagR}:
\begin{eqnarray}
	\frac{\partial S}{\partial t} + \frac{\partial J_S}{\partial x} &=& -\beta SI  \,,	\label{eq.SIR_kinetic1}		\\ 
	\frac{\partial I}{\partial t} + \frac{\partial J_I}{\partial x} &=& \beta SI  -\gamma I\,,\\ 
	\frac{\partial R}{\partial t} + \frac{\partial J_R}{\partial x} &=& \gamma I \,, \label{eq.SIR_kinetic3}\\
	\frac{\partial J_S}{\partial t} + \lambda_S^2 \frac{\partial S}{\partial x} &=& -\beta J_SI  -\frac{J_S}{\tau_S}\,, \label{eq.SIR_kinetic4}\\ 
	\frac{\partial J_I}{\partial t} + \lambda_I^2 \frac{\partial I}{\partial x} &=& \frac{\lambda_I}{\lambda_S}\beta J_SI  -\gamma J_I -\frac{J_I}{\tau_I}\,,\\ 
	\frac{\partial J_R}{\partial t} + \lambda_R^2 \frac{\partial R}{\partial x} &=& \frac{\lambda_R}{\lambda_I}\gamma J_I -\frac{J_R}{\tau_R}\,.
\label{eq.SIR_kinetic6}
\end{eqnarray}
Let us now consider the behavior of this model in the zero-relaxation limit \cite{lions1997}. To this aim, we introduce the space dependent diffusion coefficients,
\begin{equation}
D_S=\lambda_S^2 \tau_S\,,\quad D_I=\lambda_I^2 \tau_I\,,\quad D_R=\lambda_R^2 \tau_R\,,
\label{eq:diff2}
\end{equation}
which characterize the diffusive transport mechanism of susceptible, infectious and removed, respectively. Keeping the above quantities fixed while letting the relaxation times $\tau_{\mathit{S,I,R}}\to 0$ (and so the characteristic velocities $\lambda_{\mathit{S,I,R}} \to \infty$), from equations~\eqref{eq.SIR_kinetic4}--\eqref{eq.SIR_kinetic6} we obtain three proportionality relations, one for each epidemic compartment, between the flux and the spatial derivative of the corresponding density (the so-called Fick's law):
\begin{equation}
J_S = -D_S \frac{\partial S}{\partial x}\,,\quad  J_I = -D_I\frac{\partial I}{\partial x}\,,\quad J_R = -D_R\frac{\partial R}{\partial x}\,.
\label{eq:fick}
\end{equation}
Substituting \eqref{eq:fick} into equations \eqref{eq.SIR_kinetic1}--\eqref{eq.SIR_kinetic3}, we recover the following parabolic reaction-diffusion model, widely used in literature to study the spread of infectious diseases \cite{magal2019,sun2012,berestycki2021,viguerie2020}:
\begin{eqnarray}
\frac{\partial S}{\partial t} &=&  -\beta SI +\frac{\partial}{\partial x} \left({D_S}\frac{\partial S}{\partial x}\right)\,, \label{eq:diff1}\\
\frac{\partial I}{\partial t} &=&  \beta SI-\gamma I+\frac{\partial}{\partial x} \left({D_I}\frac{\partial I}{\partial x}\right)\,, \\
\frac{\partial R}{\partial t} &=&  \gamma I+\frac{\partial}{\partial x}\left({D_R}\frac{\partial R}{\partial x}\right)\,.
\label{eq:diff3}
\end{eqnarray}
This model is therefore capable to account for different regimes, ranging from hyperbolic to parabolic, according to the space dependent values $\tau_{\mathit{S,I,R}}$ and $\lambda_{\mathit{S,I,R}}$. This feature makes it particularly suitable for describing the multiscale dynamics of human beings \cite{albi2022,bertaglia2021,bertaglia2020b}.

For rigorous results on the diffusion limit of kinetic models of the type of system~\eqref{eq.SIR_kinetic_diagS}--\eqref{eq.SIR_kinetic_diagR} we refer to \cite{lions1997,salvarani2005}.
Finally, we underline that, although the model here discussed, for simplicity of presentation, is based on a simple SIR structure, the approach can be extended naturally to more realistic compartmental models designed to take into account specific features of the infectious disease of interest, such as those specifically developed  to deal with the COVID-19 pandemic \cite{gatto2020,giordano2021}.

\subsection{APNN for the hyperbolic SIR model}
The multiscale nature of the problem poses a concrete challenge to the PINN construction, and preservation of the AP property is essential in order to obtain reliable results. 
To satisfy the AP property, we follow the approach presented in Section \ref{AP-GT}, starting from the system written in macroscopic form, defined by equations~\eqref{eq.SIR_kinetic1}--\eqref{eq.SIR_kinetic6}. After multiplying both members of each equation \rev{by} the corresponding scaling parameter $\tau_{\mathit{S,I,R}}$, we can rewrite the system in the following compact form:
\begin{equation*}
    {\tau}(x)\frac{\partial {U}(x,t)}{\partial t} + {D}(x) \frac{\partial {F}({U})}{\partial x} = {G}({U})\,, \quad (x, t)\in \Omega\,, 
\end{equation*}
where 
\small
\begin{equation*}
    \label{eq:ap-vector-components}
    {U} = 
        \begin{bmatrix}
            S \\ I \\ R \\ J_S \\ J_I \\ J_R
        \end{bmatrix}
    , \
    {\tau} = 
        \begin{bmatrix}
            1 \\ 1 \\ 1 \\ \tau_S \\ \tau_I \\ \tau_R
        \end{bmatrix}
    , \
    {D}= 
        \begin{bmatrix}
            1 \\ 1 \\ 1 \\ D_S \\ D_I \\ D_R
        \end{bmatrix}, \
    {F}({U}) = \begin{bmatrix}
            J_S \\ J_I \\ J_R \\ S \\ I \\ R
    \end{bmatrix}, \
    {G}({U}) = \begin{bmatrix}
            -\beta SI \\
            \beta SI -\gamma I \\
            \gamma I \\
            -\tau_S\beta J_S I - J_S \\
            \tau_I\frac{\lambda_I}{\lambda_S}\beta J_S I - \tau_I\gamma J_I - J_I \\
            \tau_R\frac{\lambda_R}{\lambda_I}\gamma J_I - J_R
        \end{bmatrix}.
\end{equation*}
\normalsize
We consider ${U}_{\mathit{NN}}(x,t;\theta)$ to be the output vector of a DNN with inputs $x$ and $t$ and trainable parameters $\theta$, which approximate the solution of our system: ${U}(x,t) \approx {U}_{\mathit{NN}}(x,t;\theta)$. Then, we define the residual vector term
\begin{equation}
    \label{eq:pinn-ap-residue}
    \mathcal{R}^{\tau}_{\mathit{NN}} = {\tau}\frac{\partial {U}_{\mathit{NN}}}{\partial t} + {D} \frac{\partial {F}({U}_{\mathit{NN}})}{\partial x} - {G}({U}_{\mathit{NN}})
\end{equation}
and embed it into the loss function of the neural network to obtain the sought APNN~\cite{bertaglia2022a}. 
We omit for brevity the detailed analysis of the AP property. Anyway, in the limit $\tau_{\mathit{S,I,R}} \to 0$, $\lambda_{\mathit{S,I,R}}\to \infty$, under conditions \eqref{eq:diff2}, such analysis follows the same steps presented in Section \ref{AP-GT} for the prototype model, and $\mathcal{R}^{0}_{\mathit{NN}}$ results in agreement with the equilibrium system \eqref{eq:diff1}--\eqref{eq:diff3}.

If we consider both the epidemic parameters $\beta$ and $\gamma$ to be unknown, in the training process of the APNN we minimize the following AP-loss function, which recalls Eq. \eqref{eq:general-pinn-loss}:
\begin{equation}
    \label{eq:pinn-loss}
    \mathcal{L}(\theta,\beta,\gamma) = \omega_d^T \,{\mathcal L}_d(\theta,\beta,\gamma) + \omega_b^T\, {\mathcal L}_b(\theta,\beta,\gamma) + \omega_r^T \,{\mathcal L}_r(\theta,\beta,\gamma) \,, 
\end{equation}
where ${\mathcal L}_r$ enforces the minimization of the residual \eqref{eq:pinn-ap-residue}:
\begin{equation*}
\omega_r^T \, {\mathcal L}_r(\theta,\beta,\gamma) = \frac{\omega_r^T }{N_r} \sum_{n=1}^{N_r}  \left| \mathcal{R}_{\mathit{NN}}^{\tau} (x_r^n, t_r^n;\theta,\beta,\gamma) \right|^2\,.
\end{equation*}
The reader is referred to \cite{bertaglia2022a} for a detailed expansion of each term in this AP-loss function.

\subsection{APNN performance with epidemic dynamics}
Let us now examine through a numerical example the APNN performance in inferring unknown epidemic parameters, reconstructing the spatio-temporal dynamic and forecasting the spread of an infectious disease.
We design a numerical test with an initial condition that simulates the presence of two epidemic hot-spots, aligned in the spatial domain $\mathcal{D}=[0,20]$, with different numbers of infected individuals, distributed according to a Gaussian function,
$$ I(x,0)= \alpha_1\,e^{-(x-x_1)^2} + \alpha_2\,e^{-(x-x_2)^2}, $$
where $x_1=5$ and $x_2=15$ are the coordinates of the hot-spots, while 
$\alpha_1=0.01$ and $\alpha_2=0.0001$ define the different initial epidemic concentration in the two cities.
Assuming that there are no immune individuals at $t=0$ and that the total population is uniformly distributed in the domain, we impose
$$ S(x,0)=1 - I(x,0)\,, \qquad R(x,0)=0\,. $$
Initial fluxes are set to be in equilibrium, following \eqref{eq:fick}, and we consider periodic boundary conditions to allow both directions of connection for the two cities. 

As with the numerical test with the Goldstein-Taylor model, the numerical solution obtained with a second-order AP-IMEX Runge-Kutta finite volume method \cite{bertaglia2021,bertaglia2020b} is used as synthetic data for the training of the APNN when solving the inverse problem. Nevertheless, since data of fluxes $J_S, J_I, J_R$ are not accessible in real-world applications, we only enforce measurements of the densities $S, I, R$ in ${\mathcal L}_d$. 

We run two test cases for this problem: 
\begin{itemize}
\item in the first case, a hyperbolic configuration of speeds and relaxation parameters is considered, with $\lambda_{\mathit{S,I,R}}^2=1$ and $\tau_{\mathit{S,I,R}}=1$; 
\item in the second case, we simulate a parabolic regime setting $\lambda_{\mathit{S,I,R}}^2=10^3$ and $\tau_{\mathit{S,I,R}}=10^{-3}$. 
\end{itemize}
For each scenario, both an inverse and a forward problem are evaluated. 
PDEs residuals are enforced on $N_r=40000$ and $N_r=23600$ in the spatio-temporal domain respectively when simulating the parabolic and the hyperbolic regime, always using~20\% of training and residual points for validation purposes only. To strictly impose also the periodic boundary conditions in the neural network (accounted in ${\mathcal L}_b$ in the loss function), we employ the periodic mapping technique taken from \cite{zhang2020} in the input layer:
\begin{equation}
    \label{eq:pinn-periodic-mapping}
    { U}_{\mathit{NN}}(x, t) = { U}_{\mathit{NN}}\left(\cos(\alpha x), \sin(\alpha x), t\right)\,,
\end{equation}
where $\alpha$ is a hyperparameter controlling the frequency of the solution, here chosen equal to 3.  

To solve the problem, we adopt again a feed-forward neural network, this time with depth $8$ and width $32$, choosing \texttt{tanh} as activation function and fixing a learning rate LR=$10^{-2}$. Finally, the Adam method \cite{Adam} is employed for the optimization process. Weights associated to each loss term are listed in Table \ref{Table:SIRhyperparameters}.

\begin{table}[!b]
\caption{Weights used in the loss function of the APNN for the SIR test cases.}
    \label{Table:SIRhyperparameters}
    \centering
    \begin{tabular}{p{2.0cm}p{2.2cm}p{2.2cm}p{4.5cm}}\hline\noalign{\smallskip}
    Test regime &$(\omega_d^S,\,\omega_d^I,\,\omega_d^R)$ &$(\omega_b^S,\,\omega_b^I,\,\omega_b^R)$ &$(\omega_r^S,\,\omega_r^I,\,\omega_r^R,\,\omega_r^{J_S},\,\omega_r^{J_I},\,\omega_r^{J_R})$ \\
    \noalign{\smallskip}\hline\noalign{\smallskip} 
    Hyperbolic &(1, 100, 10) & (1, 10, 1) & (1, 100, 10, 1, 100, 10) \\
    Parabolic &(1, 100, 10) & (1, 10, 1) & (1, 10, 1, 1, 10, 1) \\
    \noalign{\smallskip}\hline\noalign{\smallskip}
    \end{tabular}
\end{table}
\begin{figure}[!t]
    \centering   
            \includegraphics[width=0.45\textwidth]{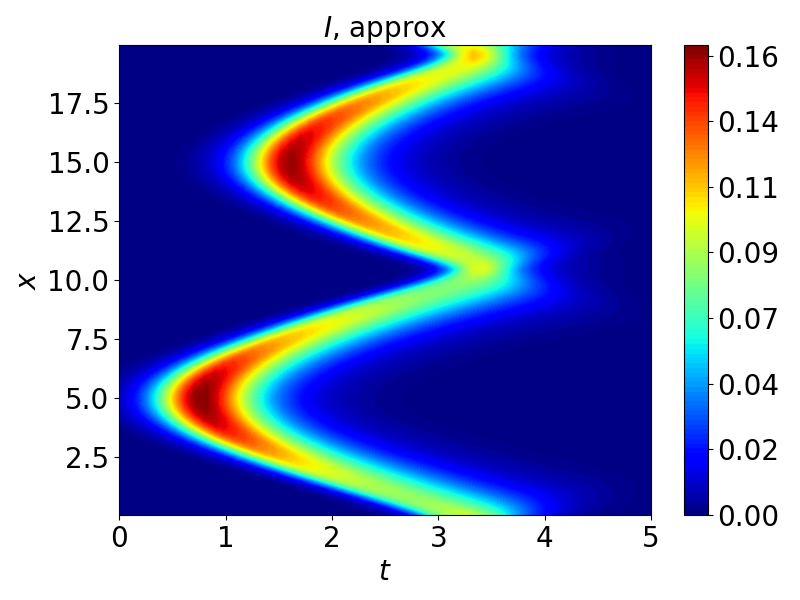}
    \includegraphics[width=0.45\textwidth]{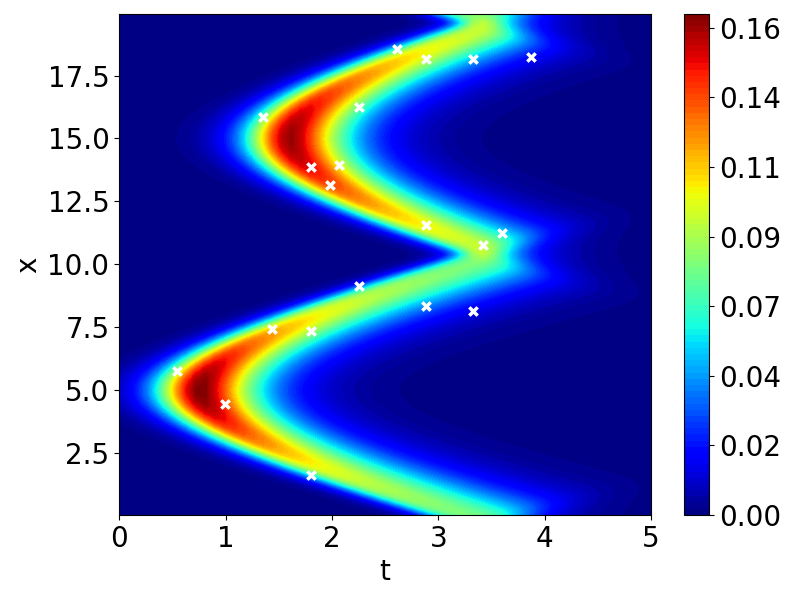}
    \includegraphics[width=0.45\textwidth]{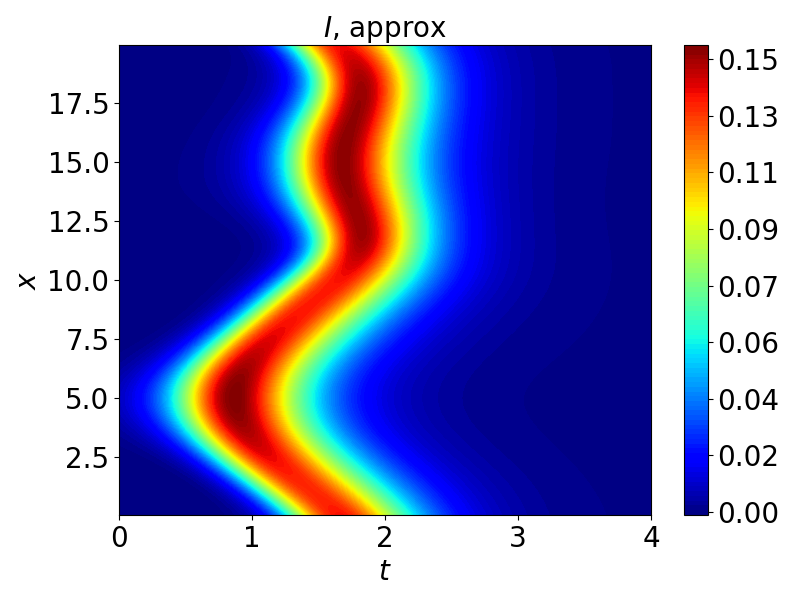}
    \includegraphics[width=0.45\textwidth]{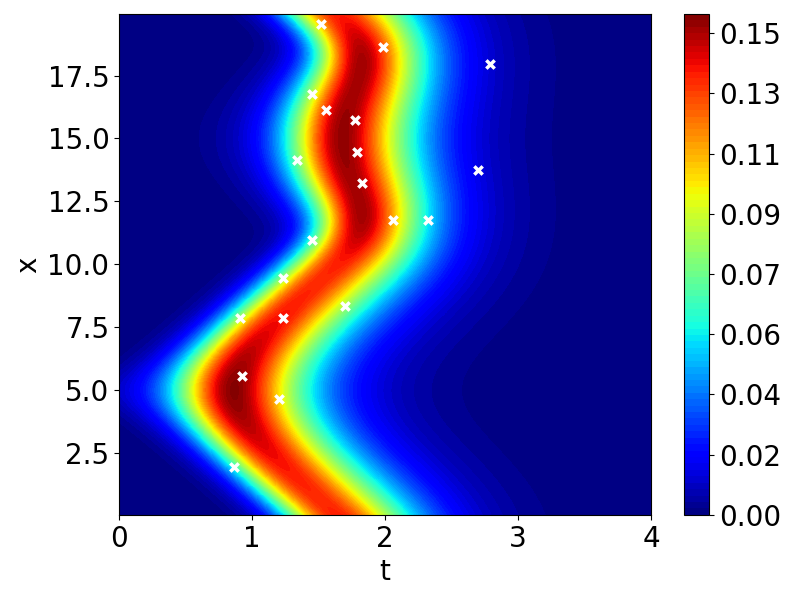}  
    \caption{Inverse problem with the SIR transport model considering a partially observed dynamic (sparse density data only and no initial conditions available) in hyperbolic regime (top) or parabolic regime (bottom). Approximation obtained with the APNN (left) and ground truth of the densities of infected $I$ (right) with the identification of the sparse data samples used for the training ($N_d=20$) marked with white crosses. }
    \label{fig:test3_inverse}
\end{figure}

\begin{table}[!b]
\caption{Inverse problem with the SIR transport model considering a partially observed dynamic (sparse density data only and no initial conditions available). Inferred results for transmission rate~$\beta$ and recovery rate $\gamma$ from a sparse measurement dataset of $N_d=20$ samples, and the relative error with respect to the ground truth values.}
\label{Table:test3_inference}
	\centering
    \begin{tabular}{p{1.7cm}p{1.5cm}p{1.9cm}p{1.7cm}p{2.4cm}p{1.8cm}}\hline\noalign{\smallskip}
    Test regime &Parameter & Ground Truth & Initial Guess & APNN Estimation & Relative Error \\
    \noalign{\smallskip}\hline\noalign{\smallskip}
    \multirow{2}{*}{Hyperbolic}&$\beta$ & 12 & 8 & 12.0126 & $1.05\times 10^{-3}$ \\
    &$\gamma$ & 6 & 3 & 6.0447 & $7.45\times 10^{-3}$ \\
    \noalign{\smallskip}\hline\noalign{\smallskip}
    \multirow{2}{*}{Parabolic}&$\beta$ & 12 & 8 & 11.9428 & $4.76\times 10^{-3}$ \\
    &$\gamma$ & 6 & 3 & 5.9772 & $3.80\times 10^{-3}$ \\
    \noalign{\smallskip}\hline\noalign{\smallskip}
    \end{tabular}   
\end{table}

\paragraph{Inverse Problems}
We initially suppose that the epidemic parameters $\beta=12$ and $\gamma = 6$ are unknown, so that the APNN is used to infer both of them as well as to reconstruct the spatio-temporal solution. To mimic the availability of data close to reality, we use a sparse dataset of only $N_d = 20$ measurements of densities $S,I,R$ for the training process, sampling the spatio-temporal points from the available dataset with probability proportional to the magnitude of $I$. Indeed, in real epidemic scenarios, data on the evolution of the infectious disease are only available in regions where the virus has already started to spread. Specifically, the probability of each spatio-temporal location $(x, t)$ chosen for the training dataset is given by
\begin{equation*}
    p(x, t) = \frac{I(x, t)}{\int_{\Omega}I(x, t)}.
\end{equation*}
Moreover, in this inverse problem, both initial and boundary conditions are regarded as unknown.

Results of the parameter inference task based on the sparse measurement dataset are reported in Table \ref{Table:test3_inference} for the test in both regimes. An excellent estimation of parameters $\beta$ and $\gamma$ can be observed with respect to the ground truth. 
At the same time, the APNN is also capable of reconstructing the correct dynamic of the epidemic spread in the whole domain besides the consistent sparsity and incompleteness of data, as shown in Figure~\ref{fig:test3_inverse} for both the test cases.

\begin{figure}[!t]
    \centering   
    \includegraphics[width=0.45\textwidth]{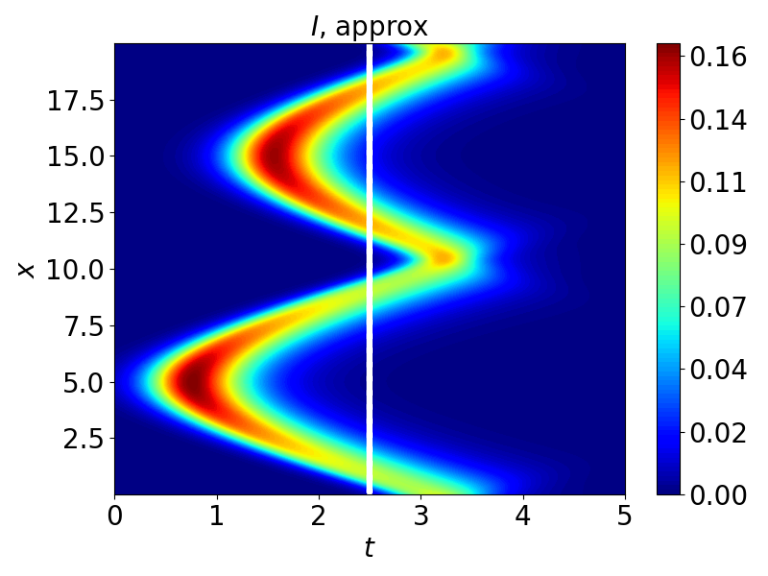}
    \includegraphics[width=0.45\textwidth]{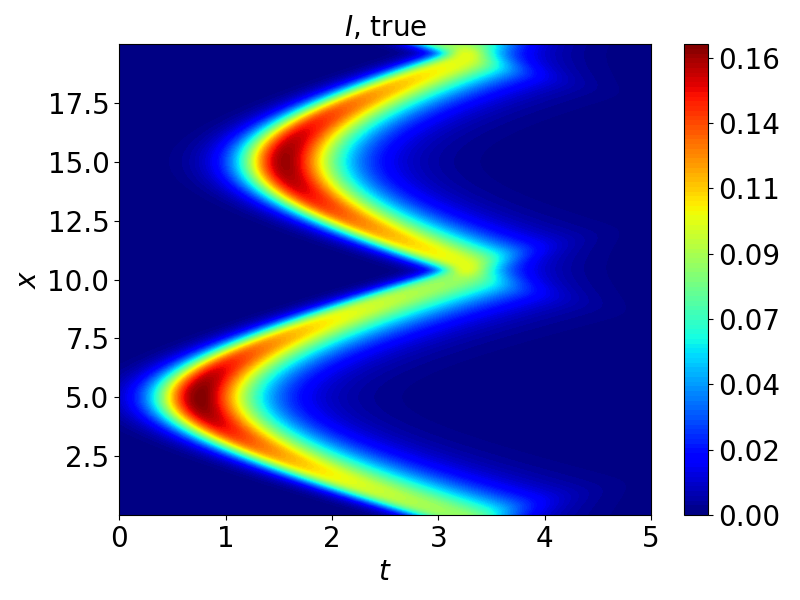}
    \includegraphics[width=0.45\textwidth]{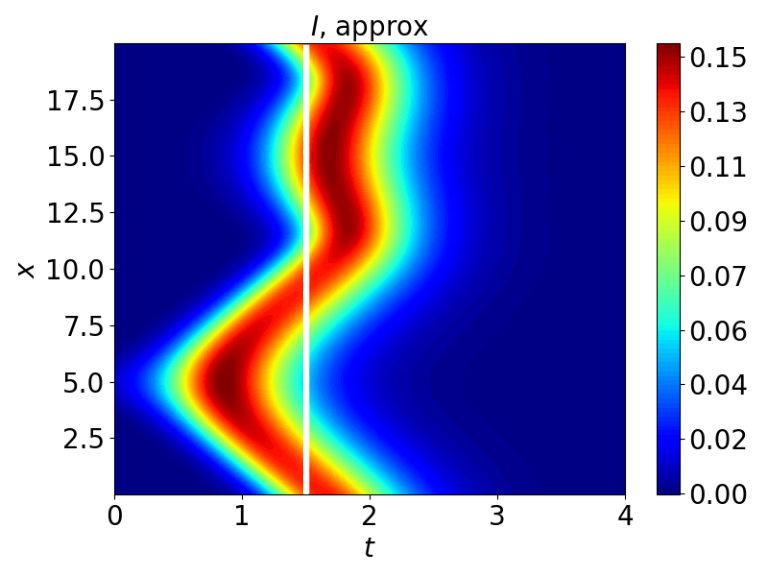}
    \includegraphics[width=0.45\textwidth]{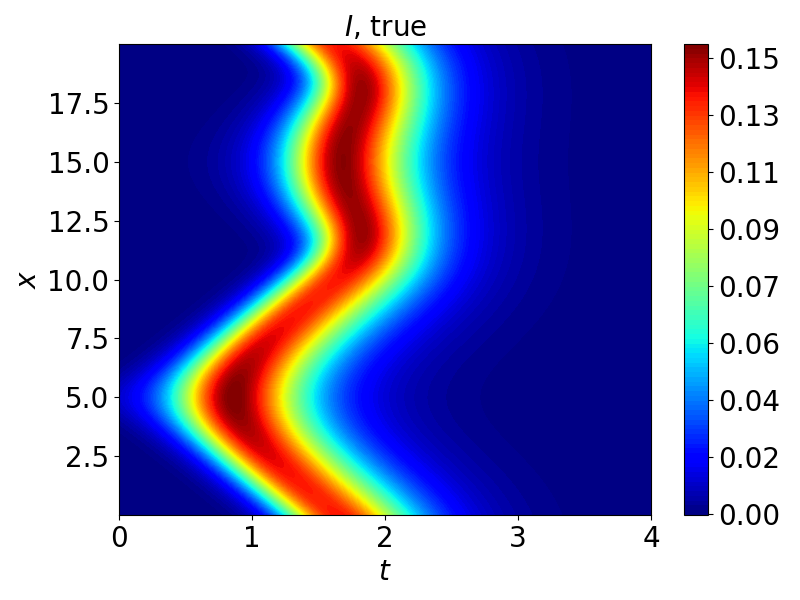}
    \caption{Forward problem with the SIR transport model considering a partially observed dynamic (sparse density data only and no initial conditions available) in hyperbolic regime (top) or parabolic regime (bottom). Approximation and forecast with measurements on a short time, with $t \in [0,2.5]$ in the hyperbolic test and $t \in [0,1.5]$ in the parabolic case (left), and ground truth of the densities of infected $I$ (right). \rev{The continuous white line on the left plots identifies the end of the data training domain.}}
    \label{fig:test3_forward}
\end{figure}

\paragraph{Forward Problems}
As a second task, we seek to investigate the forecasting capability of the APNN. In contrast to the previous inference test, in which we sampled the available measurements over the entire spatio-temporal domain, in this test we generate a training dataset over a shorter time domain and evaluate not only the correctness of the APNN approximations in this restricted time period, but also the prediction performance in subsequent time steps. In particular, 
\begin{itemize}
\item for the test in hyperbolic regime, we train the APNN with $N_d = 8500$ measurements of densities $S$, $I$ and $R$ generated from $t \in [0,2.5]$, and then evaluate the APNN performance over the complete time domain $t \in [0,5]$, with predictions in $t \in [2.5,5]$;
\item for the test in parabolic regime, we consider a training dataset of $N_d=5300$ points for densities $S$, $I$ and $R$ uniformly distributed in $t \in [0,1.5]$ and we assess the correctness of APNN approximations in $t \in [0,1.5]$ and the forecasting per- formance in $t \in [1.5,4]$.
\end{itemize}

As shown in Figure \ref{fig:test3_forward}, the approximated and predicted dynamics perfectly match the ground truth in the entire domain in both the test cases.
These results further emphasize the great capabilities of the APNN, which is also able to predict the spread of an infectious disease in multiscale regimes due to the physical knowledge of the PDEs system that has been incorporated into the loss function of the neural network in such a way as to ensure the preservation of the AP property.

Nevertheless, it is worth to remark that these predictions could only be plausible on the assumption that, at least in the short term analysis, there are no significant changes in population mobility, in the contact rate between individuals and in the contagiousness of the virus. In other words, assuming that there are no restrictions or removal of restrictions by governments and no viral variants with markedly different characteristics \rev{from} the virus already in the field begin to spread.

\section{Conclusion}
\label{sec:conclusion}
In this Chapter, we discussed a new class of Physics-Informed Neural Networks~(PINNs) that are able to appropriately deal with problems characterized by multiscale dynamics that may have several orders of magnitude difference, in particular hyperbolic problems that exhibit diffusive scaling. If we consider the residual term relating to the physical knowledge of the phenomenon under study inserted in the loss function of the PINN optimization process as a numerical approximation of the original equations, we can say that these neural networks satisfy the asymptotic-preservation (AP) property, and are therefore called Asymptotic-Preserving Neural Networks (APNNs) \cite{bertaglia2022a,jin2022}. 

Several numerical tests were presented to illustrate the performance of this new class of neural networks, addressing the problems of inverse learning and dynamics prediction. We have shown, both analytically and numerically, how APNNs provide considerably better results with respect to the different scales of the problem when compared with standard PINNs using the Goldstein-Taylor model as prototype system of equations. 

One of the application areas of particular interest, also used in this work as an example case study, concerns the investigation of the spread of infectious diseases. Indeed, it is well known that epidemic models require a delicate calibration phase of the parameters involved for their validation and subsequent use for forecasting purposes. It is perhaps less obvious to think that even the dynamic of a population, and the associated propagation of a virus, is actually defined by a system of multiscale hyperbolic equations, with individuals moving in suburban areas following a convective-transport mechanism and individuals interacting in high-density urban areas with a purely diffusive dynamic. The results obtained confirm the great usefulness and effectiveness of APNNs, even and especially when analyzing scenarios in which only little and scattered information is available. 

\section*{Acknowledgments}
The author aknowledges the support by INdAM-GNCS and by MIUR-PRIN Project 2017, No. 2017KKJP4X, \textit{Innovative numerical methods for evolutionary partial differential equations and applications}.

\bibliographystyle{unsrt_abbrv_custom}
\bibliography{APNN-hyp_postrev}
\end{document}